\documentstyle[fullpage,11pt]{amsart}

\catcode`\@=11

\long\def\@savemarbox#1#2{\global\setbox#1\vtop{\hsize\marginparwidth 
  \@parboxrestore\tiny\raggedright #2}}
\newcommand\lref[1]{\ref{#1}%
\@ifundefined{r@DisplaY #1}{}{ (#1)}}

\newcommand\fakelabel[2]{\@bsphack\if@filesw {\let\thepage\relax
   \newcommand\protect{\noexpand\noexpand\noexpand}%
\xdef\@gtempa{\write\@auxout{\string
      \newlabel{#1}{{#2}{\thepage}}}}}\@gtempa
   \if@nobreak \ifvmode\nobreak\fi\fi\fi\@esphack}

\catcode`\@=12

\def\Empty{}
\newcommand\oplabel[1]{
  \def\OpArg{#1} \ifx \OpArg\Empty {} \else
        \label{#1}
  \fi}
\newtheorem{theoremSt}{Theorem}[section]

\newtheorem{exampleSt}[theoremSt]{Example}
\newtheorem{exerciseSt}[theoremSt]{Exercise}

\newcommand\MakeStEnv[1]{
  \newenvironment{#1}[1]{
  \begin{#1St} \oplabel{##1}%
  \global\def\CrntSt{\thetheoremSt}%
}{ 
  \end{#1St} }
  \newenvironment{#1+}[1]{
  \begin{#1St} \label{##1}%
  \label{DisplaY ##1}%
  \global\def\CrntSt{\thetheoremSt}%
  \def\Labl{##1}\ifx\Labl\Empty{} \else {\em (\Labl)\,}\fi%
}{ 
  \end{#1St} }
}
\MakeStEnv{theorem}
\MakeStEnv{corollary}
\MakeStEnv{proposition}
\MakeStEnv{lemma}
\MakeStEnv{definition}
\MakeStEnv{conjecture}

\long\def\realfig#1#2#3{
\begin{figure}[htbp]
\centerline{\psfig{figure=#2}}
\caption[#1]{#3}
\oplabel{#1}
\end{figure}}

\newlength{\saveu}

\newcommand{\startproof}[1]{%
\medbreak\mbox{}\noindent{\it Proof of #1:}%
}
\newcommand{\finishproof}[1]{ 
  \def\FPArg{#1}
  \ifx\FPArg\Empty
        \newcommand\FPArg{\CrntSt}  \fi
  \smallbreak\noindent\makebox[\textwidth]{\hfill\fbox{\FPArg}}
  \medbreak\noindent
}

\newcommand{\bfheading}[1]{\par\smallskip\noindent {\bf #1}}

\newcommand\BB{{\cal B}}
\newcommand\CC{{\cal C}}

\newcommand\FF{{\cal F}}
\newcommand\GG{{\cal G}}
\newcommand\HH{{\cal H}}

\newcommand\KK{{\cal K}}
\newcommand\LL{{\cal L}}
\newcommand\MM{{\cal M}}
\newcommand\NN{{\cal N}}

\newcommand\PP{{\cal P}}

\newcommand\TT{{\cal T}}

\newcommand\PMF{{\PP\kern-2pt\MM\FF}}
\newcommand\ML{{\MM\LL}}
\newcommand\PML{{\PP\kern-2pt\MM\LL}}

\newcommand\half{{\textstyle{1\over2}}}

\newcommand\Mod{\operatorname{Mod}}
\newcommand\Area{\operatorname{Area}}

\newcommand\ep{\epsilon}

\newcommand\union{\cup}
\newcommand\intersect{\cap}
\newcommand\bbR{{\mathord{\text{I\kern-2pt R}}}}        
\newcommand\bbH{{\mathord{\text{I\kern-2pt H}}}}        

\newcommand\Z{{\bold Z}}
\newcommand\R{{\bold R}}
\newcommand\Q{{\bold Q}}

\newcommand\Hyp{{\bold H}}

\newcommand\bigrightarrow[1]{\hbox to #1{\rightarrowfill}}
\newcommand\bigleftarrow[1]{\hbox to #1{\leftarrowfill}}

\newcommand\boundary{\partial}
\newcommand\semidir{\mathrel{\hbox{\vrule depth-.03ex height1.1ex\kern-0.15em$\times$}}}

\newcommand\til{\widetilde}
\newcommand\length{\operatorname{length}}

\newcommand{\diam}{\operatorname{diam}}

\numberwithin{equation}{section}

\input{psfig}
\def\SBIMSMark#1#2#3{
 \font\SBF=cmss10 at 10 true pt
 \font\SBI=cmssi10 at 10 true pt
 \setbox0=\hbox{\SBF Stony Brook IMS Preprint \##1}
 \setbox2=\hbox to \wd0{\hfil \SBI #2}
 \setbox4=\hbox to \wd0{\hfil \SBI #3}
 \setbox6=\hbox to \wd0{\hss
             \vbox{\hsize=\wd0 \parskip=0pt \baselineskip=10 true pt
                   \copy0 \break%
                   \copy2 \break%
                   \copy4 \break}}
 \dimen0=\ht6   \advance\dimen0 by \vsize \advance\dimen0 by 8 true pt
                \advance\dimen0 by -\pagetotal
 \dimen2=\hsize \advance\dimen2 by .25 true in
  \openin2=publishd.tex
  \ifeof2\setbox0=\hbox to 0pt{}
  \else 
     \setbox0=\hbox to 3.1 true in{
                \vbox to \ht6{\hsize=3 true in \parskip=0pt  \noindent  
                {\SBI Published in modified form:}\hfil\break
                \input publishd.tex 
                \vfill}}
  \fi
  \closein2
  \ht0=0pt \dp0=0pt
 \ht6=0pt \dp6=0pt
 \setbox8=\vbox to \dimen0{\vfill \hbox to \dimen2{\copy0 \hss \copy6}}
 \ht8=0pt \dp8=0pt \wd8=0pt
 \copy8
 \message{*** Stony Brook IMS Preprint #1, #2 ***}
}

\begin{document}

\title{Geometry of the complex of curves I: Hyperbolicity}
\author{Howard A. Masur}
\author{Yair N. Minsky}
\address{University of Illinois at Chicago}
\address{SUNY Stony Brook}
\date{January 2, 1998}
\thanks{The first author was partially supported by NSF grant \#DMS
9201321. The second author was partially supported by an NSF
postdoctoral fellowship and a fellowship from the Alfred P. Sloan
Foundation.}

\maketitle
\SBIMSMark{1996/11}{October 1996}{Revised version: January 1998}
%
\newcommand\carriedby{\prec}
\newcommand\subtrack{<}
\newcommand\strongly{\carriedby\carriedby}

\renewcommand\marginpar[1]{}    

\section{Introduction}
\label{intro}
In topology, geometry and complex analysis, one 
attaches a number of interesting mathematical objects
to a surface $S$. The Teichm\"uller space $\TT(S)$ is the parameter
space of conformal (or hyperbolic) structures on $S$, up to
isomorphism isotopic to the identity. The Mapping Class Group
$\Mod(S)$ is the group of auto-homeomorphisms of $S$, up to
isotopy. The geometric and group-theoretic properties of these objects
are tied to each other via the intrinsic combinatorial
topology of $S$.

In \cite{harvey:boundary}, Harvey associated to a surface $S$
a finite-dimensional simplicial complex $\CC(S)$, called
the {\em complex of curves}, which was intended to capture some
of this combinatorial structure, and in particular to encode the
asymptotic geometry of Teichm\"uller space in analogy with Tits
buildings for symmetric spaces.  The vertices of Harvey's complex
are homotopy classes of simple closed curves in $S$, and the
simplices are collections of curves that can be realized
disjointly.  This complex was then considered by Harer
\cite{harer:stability,harer:cohomdim} from a cohomological point of view,
and by Ivanov
\cite{ivanov:complexes1,ivanov:complexes2,ivanov:complexes3} with
applications to the structure of $\Mod(S)$ (in
particular a new proof of Royden's theorem).

In this paper we begin a study of the intrinsic geometry of
$\CC(S)$, which can be made into a complete geodesic
metric space in a natural way
by making each simplex a regular Euclidean simplex of sidelength
1 (see Bridson \cite{bridson:simplicial}). 
Our main result is the following:

\begin{theorem+}{Hyperbolicity}
Let $S$ be an oriented surface of finite
type. The curve complex $\CC(S)$ is a
$\delta$-hyperbolic metric space, where
$\delta$ depends on $S$.
Except when $S$ is a sphere with 3 or fewer punctures,
$\CC(S)$ has infinite diameter.
\end{theorem+}

\noindent
(See \S\ref{hyperbolic defs} for a definition of 
$\delta$-hyperbolicity.)

We remark that in a few sporadic cases our
definition of $\CC(S)$ varies slightly from the original; see \S
\ref{complex defs}.  
Note also that we can just as well consider the 
$1$-skeleton $\CC_1(S)$ 
rather than the whole complex: $\delta$-hyperbolicity
is a quasi-isometry invariant, and $\CC(S)$ is evidently
quasi-isometric to its $1$-skeleton. 

Harer showed  \cite{harer:stability}, in the non-sporadic cases, that
$\CC(S)$ is homotopy equivalent to 
a wedge of spheres of dimension greater than 1,
and in particular is simply-connected but
not contractible. It follows that $\CC(S)$ cannot be given a
$\text{CAT}(\kappa)$ metric for any $\kappa\le 0$ (see
e.g. Ballmann \cite[\S I.4]{ballmann:spaces}).
This rules out the most simple way to 
prove $\delta$-hyperbolicity by a local argument. One might still ask
if $\CC(S)$ can be embedded quasi-isometrically in a
$\text{CAT}(\kappa)$ space for $\kappa\le 0$. If $S$ has boundary
then $\CC(S)$ embeds in 
a related {\em arc complex}, whose vertices are allowed to be
arcs with endpoints on the boundary, and which Harer proved in
\cite{harer:stability} is contractible. It is thus an interesting
question whether this complex admits a $\text{CAT}(\kappa)$ metric for
$\kappa\le 0$.

\medskip

Theorem \ref{Hyperbolicity} is motivated in part by the need to understand the
extent of an important but incomplete analogy between the
geometry of the Teichm\"uller space and that of complete,
negatively curved manifolds, or more generally of
$\delta$-hyperbolic spaces. There are many senses in which this
analogy holds, and it was exploited, for example, by Bers
\cite{bers:pseudoanosov}, Kerckhoff \cite{kerckhoff:nielsen}, and
Wolpert \cite{wolpert:nielsen}.  On the other hand, Masur
\cite{masur:teichgeo} showed that (except for the simplest
cases) the Teichm\"uller metric
on $\TT(S)$ cannot be negatively curved in a local sense, and
more recently Masur-Wolf \cite{masur-wolf} showed that it is not
$\delta$-hyperbolic.  The Weil-Petersson metric on $\TT(S)$ has
negative sectional curvatures, however they are not bounded 
away from zero \cite{masur:wpmetric,wolpert:plumbing}.

The failure of $\delta$-hyperbolicity in $\TT(S)$ is closely
related to the presence of infinite diameter regions where the
metric on $\TT(S)$ is nearly a product (let us consider from now on
only the Teichm\"uller metric on $\TT(S)$). Fixing a small $\ep_0>0$, let 
$$
H_\alpha = \{ x\in \TT(S): Ext_x(\alpha) \le \ep_0\}
$$
denote the region in $\TT(S)$
where a simple closed curve $\alpha$ has small extremal length (see Section
\ref{outline} for definitions). Then
(see Minsky \cite{minsky:extremal}) the Teichm\"uller metric in
this region is approximated by a product of
infinite-diameter metric spaces, and so cannot be
$\delta$-hyperbolic.

As a consequence of the Collar Lemma (see e.g.
\cite{keen:collar,buser:surfaces}), when $\ep_0$ is sufficiently
small the intersection pattern of the family $\{H_\alpha\}$ is
exactly encoded by the complex $\CC(S)$ (it is the {\em nerve} of
this family).  Thus, one
interpretation of our main theorem is that the regions
$\{H_\alpha\}$ are the only obstructions to hyperbolicity, and
once their internal structure is ignored, the way in which they
fit together is hyperbolic. 

This can be made precise by Farb's notion of {\em relative
hyperbolicity} \cite{farb:thesis}, and in Section \ref{relative
hyperbolicity} we will prove:

\begin{theorem+}{Relative Hyperbolicity 1}
The Teichm\"uller space $\TT(S)$ is relatively hyperbolic with
respect to the family of regions $\{H_\alpha\}$.
\end{theorem+}

A similar discussion can be carried out for the mapping class
group. A group is {\em word hyperbolic} if its Cayley graph is
$\delta$-hyperbolic. It is known that 
a group acting by isometries
on a $\delta$-hyperbolic space with finite point stabilizers and
compact quotient must itself be word-hyperbolic, and it is plain
that $\Mod(S)$ acts isometrically on $\CC(S)$, with compact
quotient. Nevertheless, $\Mod(S)$ is known not to be
word-hyperbolic for all but the simplest cases, because it
contains high-rank abelian subgroups. This is not a
contradiction, because the action on $\CC(S)$ has infinite point
stabilizers. 

Indeed, abelian subgroups in $\Mod(S)$ are generated by elements
that stabilize disjoint subsurfaces, and in particular their
boundary curves, and hence are 
``invisible'' from the point of view of coarse geometry of the
complex. One can formalize this intuition as we did with
Teichm\"uller space by considering subgroups of $\Mod(S)$ fixing
certain curves, and their cosets. In Section \ref{relative
hyperbolicity} we will carry this out and prove:

\begin{theorem+}{Relative Hyperbolicity 2}
The group $\Mod(S)$ is relatively hyperbolic with respect to
left-cosets of a finite collection of stabilizers of curves.
\end{theorem+}

Farb shows in \cite{farb:thesis} that relative-hyperbolicity
results such as these are useful in converting information about
subgroups (such as automaticity) to information about the full
groups. Although his work does not apply directly to our situation,
it is nonetheless possible to use the results of this paper as the
first step in an inductive analysis of the structure of the
Mapping Class Group. Such an analysis will be carried out in
\cite{masur-minsky:complex2}.

\medskip

We remark finally that although our main theorem has essentially a
topological statement, the proof we have found uses Teichm\"uller
geometry in an essential way. It would be very interesting to find a
purely combinatorial proof. In particular, it would be nice to have
an effective estimate of the constant $\delta$, which our proof does
not provide since it depends on bounds obtained from a 
compactness argument in the Moduli space.

\section{Outline of the Proof}
\label{outline}
In this section, after describing some necessary background, we
will give an outline of the proof of the Hyperbolicity Theorem
\ref{Hyperbolicity}, which  
reduces it to a number of assertions. These assertions will then be
proven in  sections \ref{nesting section} through \ref{hyperbolic}.

\subsection{Hyperbolicity}
\label{hyperbolic defs}
A geodesic metric space $X$ is a path-connected metric space in
which any two points $x,y$ are connected by an isometric image of
an interval in the real line, called a geodesic and denoted
$[xy]$ (we use this notation although $[xy]$ is not required to
be unique).

We say that $X$ satisfies the {\em thin triangles condition} if
there exists some $\delta\ge 0$ such that, for any $x,y$ and
$z\in X$ the geodesic $[xz]$ is contained in a
$\delta$-neighborhood of $[xy]\union[yz]$. This is one of several
equivalent conditions for $X$ to be {\em $\delta$-hyperbolic} in
the sense of Gromov, or {\em negatively curved} in the sense of
Cannon. (We remark that there are formulations of hyperbolicity
that do not require $X$ to be a geodesic space, but we will not
be concerned with them here. See Cannon \cite{cannon:negative},
Gromov \cite{gromov:hypgroups} and also
\cite{bowditch:hyperbolicity,c-d-p,ghys-harpe}.)

Important examples of hyperbolic spaces are the classical
hyperbolic space $\Hyp^n$, all simplicial trees (here
$\delta=0$), and Cayley graphs of fundamental groups of closed
negatively curved manifolds.

We note also that every finite-diameter space is trivially
$\delta$-hyperbolic with $\delta$ equal to the diameter, which is the
reason we must check that the complex of curves has infinite diameter.

\subsection{The complex of curves}
\label{complex defs}
Let $S$ be a closed surface of genus $g$ with $p$ punctures. 
Except in the sporadic cases mentioned below, 
define a
complex $\CC(S)$ as follows: $k$-simplices of $\CC(S)$ are
$(k+1)$-tuples $\{\gamma_0,\gamma_1,\ldots,\gamma_k\}$ of distinct
non-trivial homotopy classes of simple, non-peripheral closed curves,
which can be realized disjointly.  This complex is obviously
finite-dimensional by an Euler characteristic argument, and is
typically locally infinite.

\subsection*{Sporadic cases}
In a number of cases $\CC(S)$ (and hence our main theorem) is
either trivial or already well-understood.  When $S$ is a sphere
($g=0$) with $p\le 3$ punctures, the complex is empty. In this
case we can say Theorem \ref{Hyperbolicity} holds vacuously.  When $g=0$ and
$p=4$, or $g=1$ and $p\le 1$, Harvey's complex has no edges, and
is just an infinite set of vertices. In these cases it is useful
to alter the definition slightly, so that edges are placed
between vertices corresponding to curves of smallest possible
intersection number (1 for the tori, 2 for the sphere). When this
is done, we obtain the familiar {\em Farey graph}, for which
Theorem \ref{Hyperbolicity} is fairly easy to prove. See
\cite{minsky:taniguchi} for an exposition of this case.

For the remainder of the paper we exclude the surfaces
with $g=0, p\le 4$ and $g=1,p\le 1$, which we call {\em sporadic}. 

In all other cases, 
the dimension of the complex is easily computed to be $3g-4+p$, which
in particular is at least 1. 
Letting $\CC_k$ denote the
$k$-skeleton of $\CC$, we focus on the graph $\CC_1$.
We turn $\CC_1$ into a metric space by specifying that
each
edge has length 1, and we denote by $d_\CC$ the distance
function obtained by taking shortest paths. Note also that $\CC_1$ is
a geodesic metric space. 

For $\alpha,\beta\in\CC_0(S)$, let $i(\alpha,\beta)$ denote the
geometric intersection number of 
$\alpha$ with $\beta$ on $S$, which is equal to the number of transverse
intersections of their geodesic representatives in a hyperbolic metric
on $S$. 

\begin{lemma}{connected} 
If $S$ is not sporadic,
$\CC_1$ is connected.  Moreover for any two curves
$\alpha,\beta$, $d_\CC(\alpha,\beta)\leq 2i(\alpha,\beta)+1$.
\end{lemma}

\bfheading{Remark.} In fact for large $d_\CC$ a better estimate is that
$i(\alpha,\beta)$ is at least exponential in $d_\CC$, as we shall see in
Section \ref{nesting section}.

\begin{pf} Assume that $\alpha$ and $\beta$ are realized
with minimal intersection number
If $i(\alpha,\beta)=1$ then a regular neighborhood of
$\alpha\union\beta$ is a punctured torus whose boundary $\gamma$ must
be nontrivial and nonperipheral since the torus and punctured torus
are excluded. Since $\gamma$ is disjoint from both $\alpha$ and $\beta$,
$d(\alpha,\beta)=2$.

For $i(\alpha,\beta)\geq 2$, fixing two points of
$\alpha\intersect\beta$ adjacent in $\alpha$ there are two distinct ways to do
surgery on these points, replacing a segment of $\beta$ with 
the segment of $\alpha$ between them,  producing
homotopically nontrivial simple curves $\beta_1,\beta_2$
such that $i(\alpha,\beta_j)\leq i(\alpha,\beta)-1$. 
If the two intersections agree in orientation then $i(\beta,\beta_j)=1$,
and neither $\beta_j$
can be peripheral (if it bounds a punctured disk then $\alpha$ must
enter it and has a non-essential intersection with $\beta$).
Thus $d(\alpha,\beta) \le 2 + d(\alpha,\beta_j)$ and we are done by
induction.

If the two intersections have opposite signs then actually
$i(\alpha,\beta_j)\leq i(\alpha,\beta)-2$, and 
$i(\beta,\beta_j)=0$ for $j=1,2$. Thus if at least one $\beta_j$ is
nonperipheral, we again apply induction (and get a better estimate
than above). If both $\beta_1,\beta_2$ are peripheral then $\beta$
must bound a twice punctured disk on the side containing the $\alpha$
segment between the intersections. Thus consider a segment of $\alpha$
between intersections, which is adjacent to this one. If it also falls
into the last category then $\beta$ bounds a twice punctured disk on
its other side too, and $S$ must be a 4-times punctured sphere, which
has been excluded.
\end{pf}

\subsection{Teichm\"uller space}
\label{teich defs}
An analytically finite conformal structure on $S$ is an identification
of $S$ with a closed Riemann surface minus a finite number of points. 
Let $\TT(S)$ denote the Teichm\"uller space of 
analytically finite conformal structures on $S$, modulo conformal
isomorphism isotopic to the identity.

Given an element $x\in\TT(S)$ and a simple closed curve $\alpha$ in
$S$, we recall that the {\em extremal length} $Ext_x(\alpha)$ 
is the reciprocal of the largest conformal modulus of an embedded
annulus in $S$ homotopic to $\alpha$.
We remark also that an alternate definition is 
$Ext_x(\alpha) = \sup_\sigma |\alpha^*|^2_\sigma$
where $\sigma$ ranges over conformal metrics  of area 1 on $(S,x)$,
and $|\alpha^*|_\sigma$ denotes $\sigma$-length of a shortest
representative of $\alpha$.
(See e.g. Ahlfors \cite{ahlfors:invariants}.)

The Teichm\"uller metric $d_\TT$ on $\TT(S)$ can be defined in terms of maps
with minimal quasiconformal dilatation, but for us it will be useful
to note Kerckhoff's characterization \cite{kerckhoff}:
\begin{equation}
\label{kerckhoff theorem}
d_\TT(x,y) = \half\log \sup_{\alpha\in\CC_0(S)} { Ext_y(\alpha)\over
Ext_x(\alpha) }.
\end{equation}

A holomorphic quadratic differential $q$ on a Riemann surface is a
tensor of the form $\varphi(z)dz^2$ in local coordinates, with
$\varphi$ holomorphic. Away from zeroes, a coordinate $\zeta$ can be
chosen so that $q = d\zeta^2$, which determines a Euclidean metric
$|d\zeta^2|$ together with a pair of orthogonal foliations parallel to
the real and imaginary axes in the $\zeta$ plane. These are
well-defined globally and are called the {\em horizontal} and {\em
vertical} foliations, respectively.  (See Gardiner \cite{gardiner} or
Strebel \cite{strebel}.) 

Geodesics in $\TT(S)$ are determined by quadratic differentials. Given
$q$ holomorphic for some $x\in\TT(S)$, for any $t\in\R$ we consider
the conformal  structures obtained by scaling the horizontal
foliation of $q$ by a factor of $e^t$, and the vertical by
$e^{-t}$. The resulting family, which we write $L_q(t)$, is a geodesic
parametrized by arclength.

For a closed curve or arc $\alpha$ in $S$, denote by
$|\alpha|_q$ its length in the $q$ metric. Let $|\alpha|_{q,h}$
and $|\alpha|_{q,v}$ denote its horizontal and vertical lengths,
respectively, by which we mean the total lengths  of the (locally
defined) projections of $\alpha$ to the horizontal and vertical
directions of $q$. 

Finally we note that the variation of
horizontal and vertical lengths is given by 
\begin{equation}
\label{variation h}
|\alpha|_{q_t,h} = |\alpha|_{q_0,h} e^t
\end{equation}
and 
\begin{equation}
\label{variation v}
|\beta|_{q_t,v} = |\beta|_{q_0,v} e^{-t}.
\end{equation}

\subsection{The proof of the Hyperbolicity Theorem.}
\label{main proof outline}
One way to prove hyperbolicity is to find a
class of paths with the following contraction property:

\begin{definition}{contraction property}
Let $X$ be a metric space. 
Say that a path $\gamma:I\to X$ (where $I\subset \R$ is some
interval,
possibly infinite)
has the {\em contraction property} if there
exists $\pi:X\to I$ and constants $a,b,c>0$ 
such that:

\begin{enumerate}
\item \label{identity on gamma}
For any $t\in I$, $\diam(\gamma([t,\pi(\gamma(t))])) \le c$
\item \label{quasi-lipschitz}
If $d(x,y) \le 1$ then $\diam\gamma([\pi(x),\pi(y)]) \le c$.
\item \label{contraction}
If $d(x,\gamma(\pi(x))) \ge a$ and
$d(x,y) \le b d(x,\gamma(\pi(x)))$, then 
$$ \diam \gamma[\pi(x),\pi(y)] \le c.$$
\end{enumerate}
(Here for $s,t\in\R$ we take
$[s,t]$ to mean the interval with endpoints $s,t$ regardless of
order.)
\end{definition}

One should think of this property in analogy with closest-point
projection to a geodesic in $\Hyp^n$. 
Condition (1) is a coarsening of the requirement that points in
$\gamma(I)$ be fixed. Condition (2) states that the projection is
coarsely Lipschitz. Condition (3) is the most important, stating that
the map is, in the large, strongly contracting for points far 
away from their images in $\gamma(I)$. Note that this holds in
$\Hyp^n$ for $b=1$. 

Note also that we give $\pi$ as a map to the parameter interval $I$
rather than its image, in order to avoid requiring anything about the
speed of parametrization of $\gamma$: for example $\gamma$ is allowed
to be constant for long intervals, and on the other hand it need not
be continuous.

We say that a family $\Gamma $ of paths has the
contraction property   if every $\gamma\in\Gamma$ has the
contraction
property, with respect to a 
uniform $a,b,c>0$.

Call a family of paths {\em coarsely transitive}
if there exists $D\ge 0$ such that for any $x$ and $y$ with
$d(x,y)\geq D$ there is a path in the family joining $x$ to $y$.
In section \ref{hyperbolic} we will prove the following theorem,
which is
probably well-known. 

\begin{theorem}{contraction implies hyperbolicity}
If a geodesic metric space $X$ has a coarsely transitive path family
$\Gamma$ with
the contraction property then $X$ is hyperbolic.
Furthermore, the paths in $\Gamma$ are uniformly quasi-geodesic.
\end{theorem}
(See \S\ref{hyperbolic} for the definition of quasi-geodesic in
this
context).

Our family of paths will be constructed using Teichm\"uller
geodesics,
in the following manner. There is a natural map
$\Phi$ from $\TT(S)$ to finite subsets of $\CC(S)$, assigning to
any
$x\in\TT(S)$ the set of curves of shortest $Ext_x$ 
(extremal length is convenient for us, though hyperbolic will do as well). 
A geodesic in $\TT(S)$ traces out, via $\Phi$, a path in $\CC(S)$
up to some bounded ambiguity. 

That is, let $q$ be a quadratic differential on  a Riemann
surface $x$
and let $L_q:\R\to \TT(S)$ be the corresponding Teichm\"uller
geodesic (parametrized by arclength). Let a map 
$$
F \equiv F_q : \R \to \CC(S)
$$
be defined by assigning to $t$ one of the curves of
$\Phi(L_q(t))$.
The actual choices will not matter, as $\Phi(x)$ has uniformly
bounded diameter:

\begin{lemma}{Phi diam bound} 
There exists $c=c(S)$ such that $\diam_\CC \Phi(x)\le c$
for all $x\in\TT(S)$.
\end{lemma} 

\begin{pf}
There exists $e_0(S)$ such that the shortest nonperipheral curve on  $(S,x)$
has extremal length at most $e_0$. Thus Lemma
\ref{short are close} below immediately bounds the distance between
any two shortest curves, by $2e_0 + 1$. 

Note in fact that 
there exists $\epsilon_0$ such that if $x$ has a curve $\alpha$
of extremal length at most $\epsilon_0$ then any curve
intersecting  $\alpha$ has extremal length greater than
$\epsilon_0$.  In this case the diameter of $\Phi(x)$ is at most
$1$.
\end{pf}

\begin{lemma}{short are close}
For $\alpha,\beta\in\CC_0(S)$, 
if $Ext_x(\alpha) \le E$ and $ Ext_x(\beta) \le E$ for some
conformal
structure $x$ on $S$, then $d_\CC(\alpha,\beta) \le 2E+1$.
\end{lemma}
\begin{pf}
It is an elementary fact (see e.g. \cite{minsky:slowmaps})
that
$Ext_x(\alpha)Ext_x(\beta) \ge i(\alpha,\beta)^2$. Thus the
assumption
of the lemma gives $i(\alpha,\beta) \le E$. Now by Lemma
\ref{connected}, $d_\CC(\alpha,\beta) \le 2E+1$.
\end{pf}

If $q$ has a closed vertical leaf then there is a collection of
(up to homotopy) disjoint vertical curves whose extremal lengths
go to 0 as $t\to\infty$. In this case choose a fixed one of these
to be the value of $F$ as $t\to\infty$, and let this also be
denoted by $F(\infty)$.  Similarly define $F(-\infty)$ if there
are horizontal curves.

The projection for $F$ will be defined using the notion of {\em
balance}.  Recalling the notation of \S\ref{teich defs},
we say that $\beta$ is balanced with respect to
$q$ if  $|\beta^*|_{q,h} = |\beta^*|_{q,v}$, where $\beta^*$ is a
$q$-geodesic representative (it may be necessary for $\beta^*$ to
go through punctures -- see \S\ref{basic qd}).

We note that $\beta^*$ is also geodesic with respect to any
$q_t$.  Since $|\cdot|_{q_t,h}$ and $|\cdot|_{q_t,v}$ vary like
$e^t$ and $e^{-t}$ (by (\ref{variation h}) and (\ref{variation v})), if
$\beta^*$ is not entirely vertical or horizontal with respect to
$q$ there is a unique $t$ for which $\beta$ is balanced, and this
is also the minimum of the quantity
$|\beta^*|_{q_t,h}+|\beta^*|_{q_t,v}$. We observe also that,
since the $q$-length of $\beta^*$ is estimated by 
$$ 
{1\over\sqrt 2}(|\beta^*|_{q,h}+|\beta^*|_{q,v}) \le |\beta^*|_q
\le |\beta^*|_{q,h}+|\beta^*|_{q,v},
$$ 
the minimum of $|\beta^*|_{q_t}$ also occurs within bounded
distance (in fact $\half\cosh^{-1}\sqrt 2$) of the balance point.
(Compare with the projection used in \cite{minsky:projections}).

Let $\CC_b=\CC_b(q)$ denote the set of simple closed curves that
are not
entirely horizontal or vertical for $q$.  We define
$\pi=\pi_q:\CC_0\to \R$
as follows: for $\beta\in \CC_b$ let $\pi(\beta)$ be
the unique $t$ for which $\beta$
is
balanced for $q_t$. For $\beta\in
\CC\setminus \CC_b$ 
let $\pi(\beta)$ be $+\infty$ if $\beta$ is vertical, and
$-\infty$ if
$\beta$ is horizontal. (As above, in this case $F(\pm\infty)$
makes sense).

Suppose now $d(\alpha,\beta)\geq 3$.  Then $\alpha$ and $\beta$
{\em fill $S$}, in that there is no $\gamma$ disjoint from both.  
There is therefore a 
quadratic differential $q$ whose nonsingular vertical leaves are
homotopic to $\alpha$ and whose nonsingular horizontal leaves are
homotopic to $\beta$ (see \cite[Expos\'e 13]{travaux}). Then
$F_q(+\infty)=\alpha$ and 
$F_q(-\infty)=\beta$. This shows that the family $\{F_q\}$ is coarsely
transitive.

Hyperbolicity will therefore be a consequence of Theorem
\ref{contraction implies hyperbolicity} and
the following:
\begin{theorem+}{Projection Theorem}
The path family $\{F_q\}$ satisfies
the contraction property with the projections $\pi_q$ defined above.
\end{theorem+}
The proof of this theorem will be given in section
\ref{projection proof}.

We will begin in Section \ref{nesting section} by developing tools
for controlling distances between curves in $\CC_0(S)$. Using Thurston's
train-track coordinates, we will analyze a covering of $\CC_0(S)$ by
a family of polyhedra which have the property that 
a point contained in a deeply nested sequence of polyhedra will be a
definite distance from any point outside the outermost polyhedron
(Lemma \ref{easy nesting}). A partial converse to this will be the
Nesting Lemma \ref{Nesting Lemma}, which given two distant curves will
allow us to construct a nested sequence of polyhedra separating them. 
We will apply these tools to prove Lemma \ref{nesting
consequence}, which relates intersection numbers to distance in
$\CC(S)$ in a way which can be directly applied in Section
\ref{projection proof}.

Proposition \ref{positive translation distance} in Section \ref{inf diam}
will establish the infinite-diameter claim in the
Hyperbolicity Theorem \ref{Hyperbolicity}.

\section{The nested train-track argument}
\label{nesting section}

\subsection{Train-tracks}
\label{train-track intro}
We refer to Penner-Harer \cite{penner-harer} for a thorough
treatment of train-tracks, recalling here some of the
terminology.  A train-track on a surface $S$ is an embedded
1-complex $\tau$ satisfying the following properties.  Each edge
(called a branch) is a smooth path with well-defined tangent
vectors at the endpoints, and at any vertex (called a switch) the
incident edges are mutually tangent. The tangent vector at the
switch pointing toward the interior of an edge can have two
possible directions, and this divides the ends of edges at the
switch into two sets, neither of which is permitted to be
empty. Call them ``incoming'' and ``outgoing''.  The valence of
each switch is at least 3, except possibly for one bivalent
switch in a closed curve component. Finally, we require that the
components of $S\setminus \tau$ have negative generalized Euler
characteristic, in this sense: For a surface $R$ whose boundary
consists of smooth arcs meeting at cusps, define $\chi'(R)$
to be the Euler characteristic $\chi(R)$ minus
$1/2$ for every outward-pointing cusp (internal angle 0), plus
$1/2$ for every inward-pointing cusp (internal angle $2\pi$).
For the train track complementary regions all cusps are outward, 
so that the condition $\chi'(R)<0$ excludes annuli,
once-punctured disks with smooth boundary, or unpunctured disks
with 0, 1 or 2
cusps at the boundary. 
We will usually consider isotopic train-tracks to be the same. 

A {\em train route} is a non-degenerate smooth path in $\tau$; in
particular it traverses a switch only by passing from incoming to
outgoing edge (or vice versa).  A {\em transverse measure} on
$\tau$ is a non-negative function $\mu$ on the branches
satisfying the switch condition: For any switch the sums of $\mu$
over incoming and outgoing branches are equal.  A closed
train-route induces the counting measure on $\tau$.

A train-track is {\em recurrent} if every branch is contained in
a closed train route, or equivalently if there is a transverse
measure which is positive on every branch. 


Fixing a reference hyperbolic metric on $S$, a {\em geodesic lamination} in
$S$ is a closed set foliated by geodesics
(see e.g. \cite{kerckhoff:nielsen,hatcher}). A geodesic lamination is
{\em measured } if it supports a measure on arcs transverse to its
leaves, which is invariant under isotopies preserving the leaves.
The space of all compactly supported measured geodesic laminations on
$S$, with suitable topology, is known as
$\ML(S)$, and we note that different choices of reference metric on
$S$ yield equivalent spaces.
A geodesic lamination $\lambda$ is {\em carried} on $\tau$
if there is a homotopy of $S$ taking $\lambda$ to a set of train
routes. In such a case $\lambda$ induces a 
transverse measure on $\tau$, which in turn uniquely determines
$\lambda$. The set of
measures on $\tau$ gives local coordinates on $\ML(S)$,
and in fact $\ML(S)$ is a manifold (homeomorphic to a Euclidean
space).

For a recurrent train-track $\tau$, let $P(\tau)$
denote the polyhedron of measures supported on
$\tau$. We will blur the distinction between $P(\tau)$ as a
subset of
$\ML(S)$, and as a subset of the space $\R_+^\BB$ of non-negative
functions on the branch set $\BB$ of $\tau$.

We note that $P(\tau)$ is preserved by scaling, so it is a cone
on a
compact polyhedron in projective space. However we will need to
consider actual measures in $P(\tau)$ rather than projective
classes.

By $int(P(\sigma))$ we will mean the set of weights on $\sigma$
which
are positive on every branch. Note that unless the dimension
$\dim(P(\sigma))$ is 
maximal, this is different from the interior of $P(\sigma)$ as a
subset of $\ML(S)$, which is empty.

We write $\sigma\subtrack\tau$ if $\sigma$ is a {\em subtrack} of
$\tau$; that is, $\sigma$ is a train track which is a subset
of $\tau$. We also say that $\tau$ is an {\em
extension} of $\sigma$ in this case. 
We write $\sigma\carriedby\tau$ if
$\sigma$ is {\em carried} on $\tau$, by which we mean that there
is
a homotopy of $S$ such that every train route on $\sigma$ is
taken to
a train route on $\tau$.
It is easy to see that
$\sigma\subtrack \tau$ is equivalent to $P(\sigma)$ being a
subface of
$P(\tau)$, and $\sigma\carriedby\tau$ is equivalent to
$P(\sigma)\subseteq P(\tau)$. 

Say that $\sigma$ {\em fills} $\tau$ if $\sigma\carriedby\tau$
and
$int(P(\sigma))\subseteq int(P(\tau))$. 
When both tracks are recurrent this is equivalent to saying
that every branch of $\tau$ is traversed by 
some branch of $\sigma$.
Similarly, a curve $\alpha$ fills $\tau$ if $\alpha\carriedby
\tau$  
and traverses every branch of $\tau$.

Call a train-track $\tau$ {\em large} if all the components of
$S\setminus \tau$ are polygons or once-punctured polygons. We
will
also say that $P(\tau)$ is 
large in such a situation.

We say that $\sigma$ is {\em maximal} if it is not a proper subtrack
of any other track. This is equivalent to saying that all
complementary regions of $\sigma$ are triangles or punctured monogons.
Also, except in the case of the punctured torus, it is equivalent to
$\dim(P(\tau)) = \dim(\ML(S))$. In any case, maximal implies large.

A {\em vertex cycle} of $\tau$ is a non-negative measure on $\tau$
which
is an extreme point of $P(\tau)$. That is, its projective class
is a
vertex of the projectivized polyhedron.
A vertex $\mu$ is always rational, i.e. $\mu(b)/\mu(b')\in\Q$ for
$b,b'\in\BB$. This is because an irrational $\mu$ can be approximated 
by a rational $\mu'$ with the same support, and then the decomposition
$\mu = \ep\mu' + (\mu-\ep\mu')$, for $\ep$ sufficiently small, 
contradicts the assumption that $\mu$ is an extreme point. Furthermore,
up to scaling, a vertex cycle can always be realized by
the counting measure on a single, simple closed curve: the alternative
is a union of nonparallel curves, and this again gives a nontrivial
decomposition of the measure.
We will always assume that a vertex is of this form.

\medskip
\marginpar{Is this right place?}

A track $\tau$ is {\em transversely recurrent} if every branch of
$\tau$ is crossed by some simple curve $\alpha$ intersecting $\tau$ 
transversely and {\em  efficiently} -- that is, so that
$\alpha\union\tau$ has no bigon 
complementary components. For us it will only be important to note the
following equivalent geometric property, which is established in Theorem
1.4.3 of \cite{penner-harer}: $\tau$ is transversely recurrent if and
only if for any $L$ (large) and $\ep$ (small) there is a complete
finite-area hyperbolic metric on $S$ in which $\tau$ can be realized
so that all edges have length at least $L$ and curvatures at most
$\ep$ (including  at the switches).

Furthermore, we record (Lemma 1.3.3 of \cite{penner-harer})
that if $\tau$ is transversely recurrent and $\tau'$ is a subtrack of
$\tau$ or is carried in $\tau$, then $\tau'$ is also transversely
recurrent.

We call a track {\em birecurrent} if it is both recurrent and
transversely recurrent.

\subsection*{Splitting}
Let $\tau$ be a generic train-track (all switches are trivalent).
A {\em splitting} move is one of the three elementary moves on a
local configuration, as shown in figure \ref{splitting moves}.
The three splits are called a left split, a collision, and a
right split, and the resulting tracks in each case are carried by
$\tau$. If we place a positive measure $\mu$ on $\tau$ and use
the labelling as in the figure, we note that a positive measure
is induced on the right split track if $\mu(a)>\mu(c)$, on the
left split track if $\mu(a)<\mu(c)$, and on the collision track
if $\mu(a)=\mu(c)$. We call $a$ a {\em winner} of the splitting
operation if $\mu(a)>\mu(c)$ (note that $d$ is then also a
winner).

\realfig{splitting moves}{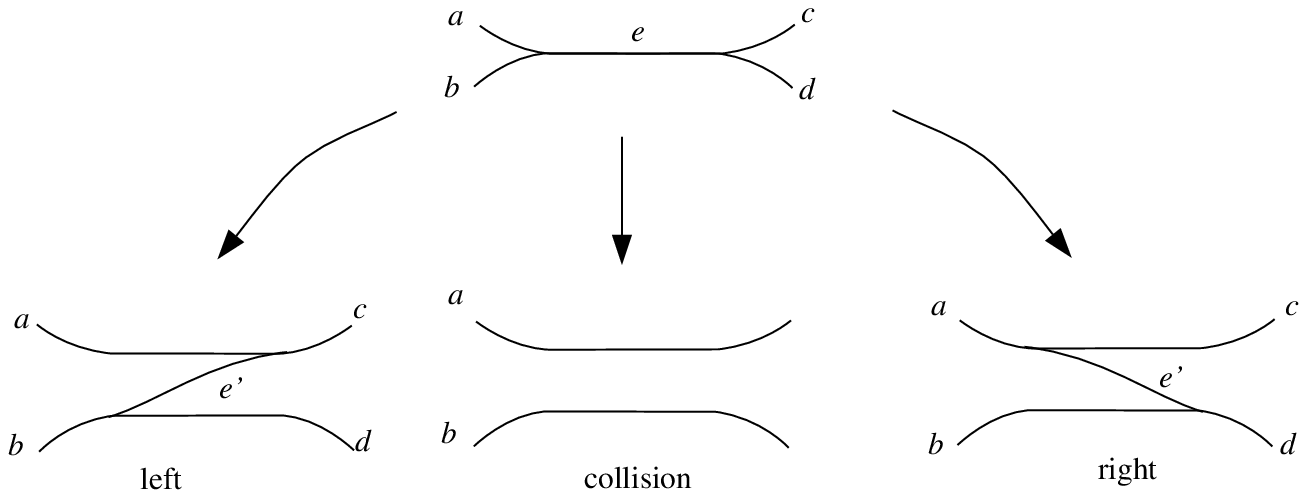}{The three ways to split
through an edge.}

Any measured lamination $\beta$ carried on $\tau$ determines a
sequence of possible splittings by the rule in the previous
paragraph, and all the resulting train-tracks carry $\beta$.
Note that if $\beta$ fills $\tau$ then it will continue to fill
the split tracks, and in particular they will all be recurrent.
This process can continue as long as the split track is not a
simple closed curve.  (We must check that for any recurrent track
that is not a simple curve there is a ``splittable'' edge, that
is one which is in the configuration of figure \ref{splitting
moves}: simply consider any transverse measure whose support is
all
of $\tau$ and take an edge of maximal weight. See
\cite{kerckhoff:simplicial}).  

When $\beta$ is a simple closed curve we will usually terminate
the
sequence as soon as we reach a track for which $\beta$ is a
vertex.

Note also that if $\sigma$ is a right or left splitting of $\tau$
then
$P(\sigma)$ and $P(\tau)$ must share at least one vertex.  To see
this, note that $P(\sigma)$ is one of the pieces obtained by
cutting
$P(\tau)$ by a hyperplane. Such a subset always contains a vertex
of
the original polyhedron.  In the case of a collision splitting,
we at
least see that $\sigma$ is a subtrack of a track that shares a
vertex
with $\tau$.

Finally, we note that when $\tau$ is
not generic, each switch can be slightly perturbed (``combing''
in
\cite{penner-harer}) to yield a generic track carrying $\tau$ 
which
carries the same set of laminations. 

\bfheading{A basic observation.}
Although it is relatively easy to understand geometrically when
pairs
of curves are a distance at most 2 in $\CC_1$, larger distances
are
more subtle to detect. This entire section is motivated by the
observation that, if
$\alpha$
and $\beta$ are disjoint curves (distance $1$ in $\CC_1$) and
$\alpha$ is carried on a
maximal
train-track $\sigma$ in such a way that it passes through every
branch, then $\beta$ is also carried on $\sigma$.  
In other words 
\begin{equation}\label{nesting 1}
\NN_1(int(P(\sigma))) \subset P(\sigma),
\end{equation}
where $\NN_1$ denotes a radius 1 neighborhood in $\CC_1$. 
This
implies inductively that if $\tau_j, j=0,\ldots,n$ is a sequence
of maximal tracks such that $P(\tau_j)\subset int( P(\tau_{j-1}))$
and $\beta_j,j\geq 1$ is a sequence in $\CC_1$ 
such that $\beta_j$ is in $int (P(\tau_{j-1}))$ but not 
in $P(\tau_j)$, then
$$d_{\CC_1}(\beta_1,\beta_j) \geq j. $$  The issue is more
subtle when the $\tau_j$ are not maximal, or equivalently if
they are maximal but the $\beta_j$ are carried on a face. 
This leads to a discussion of diagonal extensions and Lemma \ref{1-nbd
  is closure}
which generalizes the above inequality.  A partial converse is
given in Lemma \ref{Nesting Lemma}.  These two Lemmas are then
applied to prove Lemma \ref{nesting consequence}.

\subsection*{Diagonal extensions}
Let $\sigma$ be a large track. A {\em diagonal extension} of
$\sigma$ is a
track $\kappa$ such that $\sigma\subtrack\kappa$ and every branch
of
$\kappa\setminus \sigma$ is a {\em diagonal} of $\sigma$: that
is, its
endpoints terminate in corner of a complementary region of
$\sigma$. 
It is easy to see that if $\sigma$ is transversely recurrent then so
is any diagonal extension -- after realizing $\sigma$ with long edges 
with nearly zero curvature, its complementary regions are nearly
convex and we can make the diagonals nearly geodesic too.

Let $E(\sigma)$ denote the set of all recurrent
diagonal extensions of $\sigma$. Note
that it is a finite set, and
let $PE(\sigma)$ denote $\bigcup_{\kappa\in E(\sigma)}
P(\kappa)$.

Further, let us define $N(\tau)$ to be the union of $E(\sigma)$
over
all large recurrent subtracks $\sigma\subtrack\tau$. Define $PN(\tau) =
\bigcup_{\kappa\in N(\tau)} P(\kappa)$.  In some sense this
should be
thought of as a ``neighborhood'' of $P(\tau)$; compare Lemma
\ref{1-nbd is closure} and (\ref{1-nbd for PN}).

Let $int(PE(\sigma))$ denote the set of
measures $\mu\in PE(\sigma)$
which 
are positive on every branch of $\sigma$. 
We also define
$int (PN(\tau)) = \bigcup_\kappa int (PE(\kappa))$, where
$\kappa$
varies
over the large recurrent subtracks of $\tau$.

If a lamination $\mu\in \ML(S)$ has complementary regions which are
ideal polygons or once-punctured ideal polygons, then an
$\ep$-neighborhood of $\mu$, for $\ep$ sufficiently small, gives rise
to a train track $\sigma$ (see 
Penner-Harer \cite{penner-harer}) which must be large, and $\mu$ puts
positive weight on every branch of $\sigma$.
One of the
properties of the measure topology on $\ML(S)$ is then that in a small
enough neighborhood of $\mu$ all laminations are carried on some
diagonal extension of $\sigma$, and in particular in $int(PE(\sigma))$.
A more quantitative form of this idea is 
the following sufficient condition for containment
in $int(PE(\sigma))$.

\begin{lemma}{interior of extension}
There exists $\delta>0$ (depending only on $S$) for which the
following holds.  Let $\sigma\subtrack\tau$ where $\sigma$ is a
large
track.  If $\mu\in P(\tau)$ and, for every
branch $b$ of $\tau\setminus\sigma$ and $b'$ of $\sigma$, $\mu(b)<
\delta \mu(b')$, then $\sigma$ is recurrent and $\mu\in
int(PE(\sigma))$.
\end{lemma}

\begin{pf}
Whenever there are branches of $\tau\setminus\sigma$
which meet an edge $e$ of the boundary of a complementary region
of
$\sigma$ at other than a corner point, 
the situation can be simplified by either a {\em split} or a {\em shift}.
First, 
there may be a splitting move involving these branches which replaces
$\sigma$ by an equivalent track (also called $\sigma$). Such a move
either reduces
the number of edges of $\tau\setminus\sigma$ incident to
$\sigma$, or
moves one of them closer to a corner (see figure 
\ref{make it diagonal}).  These are the only two possibilities. 

\realfig{make it diagonal}{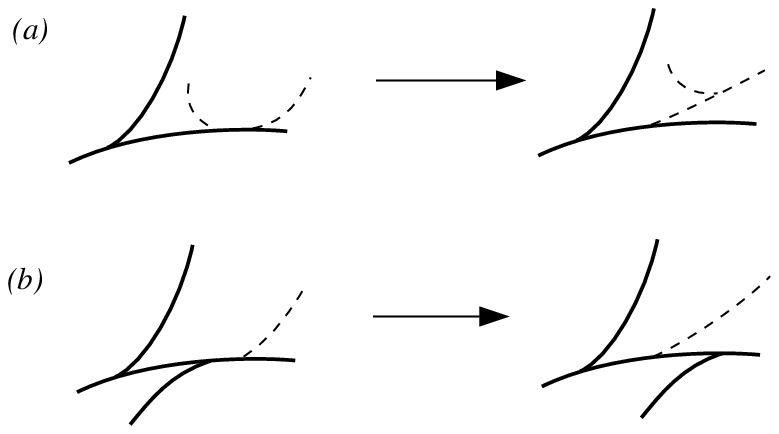}{Solid edges are in
$\sigma$, dotted edges in $\tau\setminus\sigma$. In (a), the
splitting reduces the number of dotted edges incident to
$\sigma$.  In (b), the splitting brings a dotted edge closer to a
corner.}


If a branch of $\tau\setminus\sigma$ is facing ``toward'' a corner of
a complementary region, and 
separated from it
only by non-splittable edges, we can
perform a {\em shift move} (see figure \ref{shift move}) which
takes
this branch to the corner without affecting the set of measures
carried on the track. 
Thus any sequence of such
splitting and shifting
moves must terminate after a bounded number of steps in 
a track $\tau'\in E(\sigma)$.

\realfig{shift move}{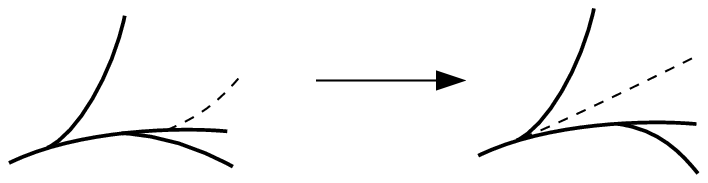}{Shifting an edge of
$\tau\setminus\sigma$ into a corner}

The measure $\mu$ will be carried on such a $\tau'$ provided it is
consistent with the sequence of moves. 
That is, whenever a
splitting is determined by a comparison between a branch $b$ of
$\tau\setminus\sigma$ and a branch $c$ of $\sigma$, the branch of
$\sigma$ must win (that is, $\mu(b) < \mu(c)$). 
After
such a splitting, there is a branch of $\sigma$ with measure
$\mu(c)-\mu(b)$.  Thus, if $k$ is the bound on the number of
splittings that take place, a sufficient condition that all the
splittings are consistent with $\mu$ is $\min_\sigma\mu >
(k+1)\max_{\tau\setminus\sigma} \mu$.

Therefore, setting $\delta = 1/(k+1)$, we  guarantee that
$\mu\in PE(\sigma)$, and furthermore that $\mu$ puts positive
measure
on each branch of $\sigma$ after the splitting, so that $\mu\in
int(PE(\sigma))$. 

Finally, we must show that $\sigma$ is recurrent. It is an easy fact
of linear algebra that, 
\marginpar{EXPLAIN BETTER}
if there is an assignment of positive weights
on the branches of $\sigma$ such that at every switch the difference
of incoming and outgoing weights is less than a fixed constant
$\delta_1$ times the minimum weight (where $\delta_1$ depends just on
the surface $S$), then these weights can be perturbed to positive
weights satisfying the switch conditions, and hence $\sigma$ is
recurrent.  It follows that, with sufficiently small $\delta$, the
measure $\mu$ restricted to $\sigma$ gives such a set of weights.
\end{pf}

The next two lemmas show that, when tracks are nested, their diagonal
extensions are nested in a suitable sense, and the way in which the
diagonal branches cover each other is controlled.

\begin{lemma}{easy nesting}
Let $\sigma$ and $\tau$ be large recurrent tracks, and suppose
$\sigma\carriedby\tau$. 
If  $\sigma$ fills $\tau$, then 
$PE(\sigma) \subseteq PE(\tau)$.
Even if $\sigma$ does not fill $\tau$, we have
$PN(\sigma)\subseteq PN(\tau)$.
\end{lemma}

\marginpar{Howie, you had some concerns here about ties emanating from
  corners, which don't seem to me to be a problem.}

\begin{pf}
  We may thicken $\tau$ slightly to get a regular neighborhood
  $\tau_\ep$, which can be foliated by short arcs called ``ties''
  transverse to $\tau$.
Then  
$\sigma$ can be embedded
in $\tau_\ep$ so that it is transverse to the ties. 

The assumption that $\sigma$ fills $\tau$ 
implies that every edge of $\tau$ is
traversed by some edge of $\sigma$, and thus $\sigma$ crosses
every
tie.  
Any component $D$ of $S\setminus\sigma$, which is a polygon or a
once-punctured polygon, must have some subset $F$ foliated by
ties.  $F$ consists of a neighborhood of the boundary and bands
joining different boundary edges. Each component of $D\setminus
F$ is isotopic to a component of $S\setminus\tau$, and the
quotient of $D$ obtained by identifying each tie to a point can
be identified with some union of complementary regions of $\tau$,
joined by train routes in $\tau$. Any diagonal edge $e$ in $D$
joining two corners of $D$ may therefore be put in minimal
position with respect to the ties (so that $e$ meets ties
transversely, and no disks are bounded by a segment of $e$ and a
tie), and hence gives rise to a train route through the union of
$\tau$ with some diagonal edges.

It follows that any diagonal extension of $\sigma$ can be carried
by a diagonal extension of $\tau$, and hence $PE(\sigma)\subset
PE(\tau)$.

Now in general, 
if $\kappa$ is a large recurrent subtrack of $\sigma$, let $\rho$ be the
smallest subtrack of $\tau$ carrying $\kappa$. Note that $\rho$
is
necessarily large and recurrent, and $\kappa$ fills $\rho$. Thus the
same argument applies to the faces, and we can conclude
$PN(\sigma)\subseteq PN(\tau)$
\end{pf}

\begin{lemma}{diagonals nesting}
Let $\sigma\carriedby\tau$ where $\sigma$ is a large recurrent
track, and let
$\sigma'\in E(\sigma)$, $\tau'\in E(\tau)$ such that
$\sigma'\carriedby\tau'$. Then any branch $b$ of
$\tau'\setminus\tau$
is traversed with bounded degree $m_0$ by branches of $\sigma'$.
The number $m_0$ depends only on $S$. 
\end{lemma}

\begin{pf}
It will suffice to show that no branch of $\sigma'$ passes
through a
branch of $\tau'\setminus\tau$ more than twice; we then obtain
$m_0$ from a topological bound on the number of branches of any
train
track in $S$.

As in the previous lemma, let $\tau'_\ep$ be a regular
neighborhood of
$\tau'$ foliated by ties 
and isotope $\sigma'$ so that it is contained in $\tau'_\ep$ and
is
transverse to the ties. Because each component of $S\setminus
\tau'$
is either a polygon with $d\ge 3$ corners or a once-punctured
polygon
with $p\ge 1$ corners, we may extend the ties to a foliation
$\FF$ of
$S$ with one index $1-d/2 \le -1/2$ singularity in each $d$-gon
and
an index $1 - p/2 \le 1/2$ singularity at the puncture of each
punctured $p$-gon. (See figure \ref{foliate complement}.)

\realfig{foliate complement}{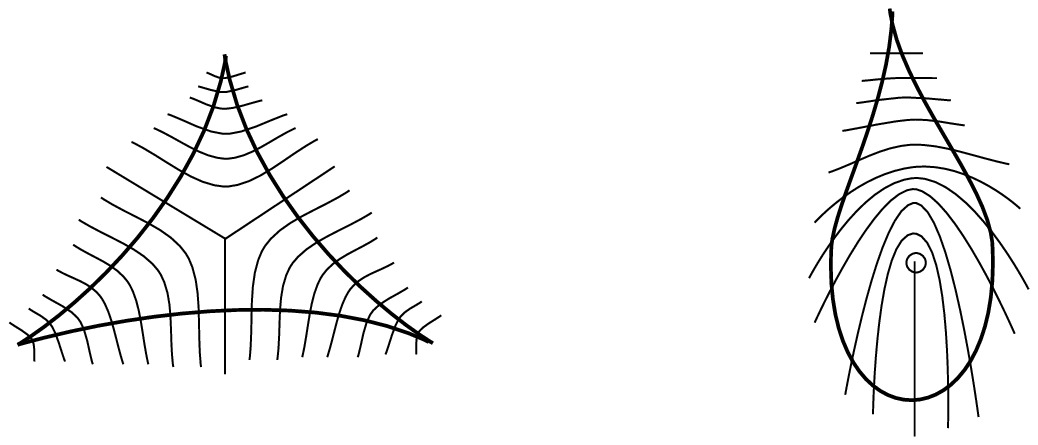}{The foliation $\FF$ in
a triangular component, and in a punctured monogon.}

Fix a tie $t$ that crosses the regular neighborhood of a branch
$b$ of
$\tau'\setminus\tau$. Since $\sigma\carriedby\tau$, no branch of
$\sigma$ crosses $t$. Suppose a branch $e$ of
$\sigma'\setminus\sigma$
crosses $t$ twice and let $t'$ be a segment in $t$ between two
successive crossings. Then $t'$ and an interval $e'\subset e$
form an
embedded loop in a complementary region $U$ of $\sigma$, which is
either a disk or once-punctured disk since $\sigma$ is large. 

Since the foliation is transverse to $e$ and parallel to $t$, we
can
see that the index of the foliation around $t'\union e'$ is
therefore
$+1$ if the ends of $e'$ meet $t$ on opposite sides, and $+1/2$
if they meet $t$ on the same side. The first case cannot occur
since the
disk bounded by $t'\union e'$ can contain at most one singularity
of
$\FF$, of index at most $1/2$. The second case can occur if
$t'\union
e'$ surrounds a puncture. 
\marginpar{Figure??}
However, in this case we can see that there
are no other points of $t\intersect e$. For if there were we
could
combine two such loops to find a loop in $U$ with index $1$,
again a
contradiction.
\end{pf}

\subsection{Nesting and $\CC$-distance}
Equation (\ref{nesting 1}) is a special
case of the following more general fact:

\begin{lemma}{1-nbd is closure}
If $\sigma$ is a large birecurrent train-track and $\alpha \in
int(PE(\sigma))$ then
$$ d_\CC(\alpha,\beta) \le 1 \implies \beta \in PE(\sigma).$$
In other words,
$$ \NN_1(int(PE(\sigma))) \subset PE(\sigma),$$
where $\NN_1$ denotes a radius 1 neighborhood in $\CC_1$. 
\end{lemma}

\begin{pf}
If $\beta\notin PE(\sigma)$ then, by the more
quantitative Lemma \ref{quantitative intersection} below,
$i(\alpha,\beta)$ is no less than the minimum weight $\alpha$
puts on any branch of $\sigma$. This is positive since $\alpha\in
int(PE(\sigma))$, but then it follows $d_\CC(\alpha,\beta) >1$.
\end{pf}

We note the following immediate consequence:
\begin{equation}
\label{1-nbd for PN}
\NN_1(int(PN(\sigma)) \subset PN(\sigma)
\end{equation}
which is 
obtained by applying Lemma \ref{1-nbd is closure} to the large
subtracks of $\sigma$.

It remains to prove the following lemma, 
which will also be used
at the end of this section.
\begin{lemma}{quantitative intersection}
Let $\tau$ be a large birecurrent track, let $\alpha$ be carried on
a diagonal extension 
$\tau'\in E(\tau)$,  and let $\beta$ be a curve not carried on
any diagonal 
extension of $\tau$.
Then $i(\alpha,\beta)
\ge \min_b\alpha(b)$ where the right hand side denotes the
minimum
weight $\alpha$ puts on all branches $b$ of $\tau$.
\end{lemma}


\newcommand\tbeta{\til\beta}
\newcommand\tbp{\til\beta_+}
\newcommand\tbm{\til\beta_-}
\begin{pf}
Consider first the case that $\tau' = \tau$, so that $\alpha$ is
actually carried in $\tau$. 


Since $\tau$ is transversely recurrent, for any $\ep>0$  there is a
hyperbolic metric 
on $S$ for which all train routes through $\tau$ have curvature at
most $\epsilon$. Fixing such a metric for $\ep<1$, and 
lifting $\tau$ to a train-track $\til\tau$ in the universal cover
$\Hyp^2$, we have 
that each train route $r$ of
$\til\tau$ is uniformly quasi-geodesic, and in particular has two
distinct endpoints $\boundary r  = \{r_+,r_-\}$
on the circle $\boundary\Hyp^2$ and stays in a
uniform neighborhood of the geodesic $r^*$ connecting them. 

Choose a component $\tbeta$ of the lift of $\beta$ to $\Hyp^2$,
and 
note it is also quasi-geodesic. Let $T_\beta$ be a generator of
the 
subgroup of $\pi_1(S)$ preserving $\tbeta$.
We say that an edge $\til e$ of $\til\tau$ {\em separates
$\tbeta$
consistently} if, for any train  route $r$ passing through $\til
e$,
its endpoints $r_\pm$ separate $\tbeta_\pm$. If this occurs for
some
$\til e$ then, letting $e$ be its projection to $S$, we deduce
immediately that $i(\beta,\alpha) \ge \alpha(e)$, since $\alpha$
lifts
to a collection of train routes with $\alpha(e)$ of them passing
through $\til e$, and through each translate $T_\beta^m(\til e)$.

Thus, let us now prove that, if no edge of $\til\tau$ separates
$\tbeta$ consistently, then $\tbeta$ is carried on a diagonal
extension
of $\til\tau$ (at the end we will check that this projects down
to an
extension of $\tau$).

Each train route $r$ separates $\Hyp^2$ into two open disks, which we
call halfplanes. Note that each is contained in a uniformly bounded
neighborhood of a geodesic halfplane bounded by $r^*$. For a halfplane
$H$, let $H'$ denote $\bar H \setminus \bar r$ (where the bar denotes
closure in the closed disk), i.e. the union of $H$ with an open arc on
the boundary circle.

Let $J_+$ and $J_-$ be the components of $\boundary\Hyp^2 \setminus
\boundary\tbeta$.  Let $\HH_+$ denote the union of all halfplanes $H$
(bounded by train routes) that meet the boundary entirely
in $J_+$, and define $\HH_-$ similarly.

Any halfplane $H_+$ from $\HH_+$ must be disjoint from any halfplane
$H_-$ from $\HH_-$. Suppose otherwise: their 
intersection would be an open set which is
either bounded or meets infinity only at $\boundary\til \beta$. In the
first case it must be a bigon bounded by arcs of the train route
boundaries of $H_+$ and $H_-$. A bigon has generalized Euler
characteristic (as defined in \S\ref{train-track intro}) $\chi'=0$.
But since, as is easily checked, $\chi'$ is additive for unions of
closures of complementary regions of $\til\tau$ (see also
Casson-Bleiler \cite{casson:unpub}), this contradicts the condition
$\chi'<0$ for complementary regions of the train-track. 
If the intersection is unbounded it is a ``generalized bigon'',
i.e. a region $R$ between two train routes $r_+$ and $r_-$, which are
either biinfinite and hence parallel to $\til\beta$, or are infinite
rays emanating from a common point and asymptotic at infinity. Since
they are quasigeodesics they remain a bounded distance apart.
Either
case again contradicts the Euler characteristic condition: For any
$M>0$, the intersection of $R$ with a ball of radius $M$ contains a
union $R_M$ of complementary regions with a {\em bounded} number of cusps,
only depending on the distance between $r_+$ and $r_-$. Thus
$|\chi'(R_M)|$ is bounded. On the other hand, by additivity,
$-\chi'(R_M)$ grows linearly with $M$. This is a contradiction.

We conclude that $\HH_+$ and $\HH_-$ are disjoint.
Now let $\KK = \Hyp^2\setminus(\HH_+\union\HH_-)$. This is a closed
set, nonempty since $\Hyp^2$ is connected, and by construction must be
a union of vertices, arcs and complementary regions of $\til\tau$.

We further claim that $\KK$ is connected, and in fact {\em
  $\ep'$-convex} for 
some small $\ep'<1$, meaning that any two points in $\KK$ are connected
by a path in $\KK$ whose curvature is bounded by $\ep'$ at every
point.

To see this, note that $\KK$ is the intersection of a sequence of {\em
  closed} halfplanes $\{C_n\}$, each the complement of an open
  halfplane from $\HH_\pm$. Letting $\KK_n = \intersect_{i<n}C_i$, 
we show the claims for each $K_n$. First consider the local
  structure of $\boundary \KK_n$. 

Suppose that two branches $a$ and $b$ of $\til\tau$, which are in
$\KK_n$, meet at a
switch $s$. If they are on the same side of the switch, the cusp
region between them must be contained in $\KK_n$. For otherwise there
is a halfplane outside $\KK_n$, whose boundary passes through $s$ and
separates $a$ 
from $b$ -- so one of them is in the halfplane and not in $\KK_n$.

We conclude from this that the boundary of $\KK_n$ is comprised of
train routes meeting at {\em outward pointing} cusps. It follows that
each component $X$ of $\KK_n$ is $\ep'$-convex for some $\ep'$ which
can be made arbitrarily  small by suitable choice of $\ep$.
(For example, use the fact that a small radial neighborhood of $X$
is convex, to deform into $X$ any geodesic with endpoints in $X$).

Furthermore, we can see by induction that $\KK_n$ is
connected: $\KK_2 = C_1$, and $\KK_{n+1} =
\KK_n\intersect C_n$. If $\KK_n$ is connected, we note that the train
route boundary of $C_n$ 
must intersect $\KK_n$ in a connected (possibly infinite) arc, for
otherwise we
obtain an arc in $\Hyp^2\setminus\KK_n$ with endpoints on
$\KK_n$,
which together with a piece of $\boundary\KK_n$ must bound a
region with $\chi'\ge 0$, again a contradiction.
It follows that $\KK_{n+1}$ is connected.

The intersection of a nested sequence of
closed $\ep'$-convex sets is itself $\ep'$-convex, and in particular
connected. Thus our conclusions follow for $\KK$.

We now claim that
$int(\KK)$ is disjoint from $\til\tau$, and hence is a  union of
  complementary components.
For, if $\til e$ is any branch of
$\til\tau$, 
by our assumption $\til e$ does not consistently separate
$\tbeta$,
and hence lies on a train route both of whose endpoints are in
either
$\bar J_+$ 
or $\bar J_-$. Thus, $\til e$ is on the boundary of one of the
halfplanes comprising either $\HH_+$ or $\HH_-$, and hence in
$\bar\HH_\pm$. (Note, this is the only place where we use the
assumption).

Since $\KK$ is invariant by $T_\beta$, its closure in
$\bar\Hyp^2$ contains $\boundary
\til\beta$. This together with $\ep'$-convexity implies that $\KK$
contains an infinite path $p$ of curvature bounded by $\ep'$, connecting
$\til\beta_+$ with $\til\beta_-$. Using the above description of the local
structure of $\KK$, we see that $p$ may be taken to be a union of train
routes and paths through diagonals of complementary regions of
$\til\tau$. That is, $p$ is carried by a diagonal extension of
$\til\tau$.  Both $p$ and the extension may be assumed invariant by $T_\beta$, 
so $p$ projects to a closed
curve in $S$ homotopic to $\beta$, and carried on $\tau$ together with
a number of diagonal branches.

\realfig{Region K}{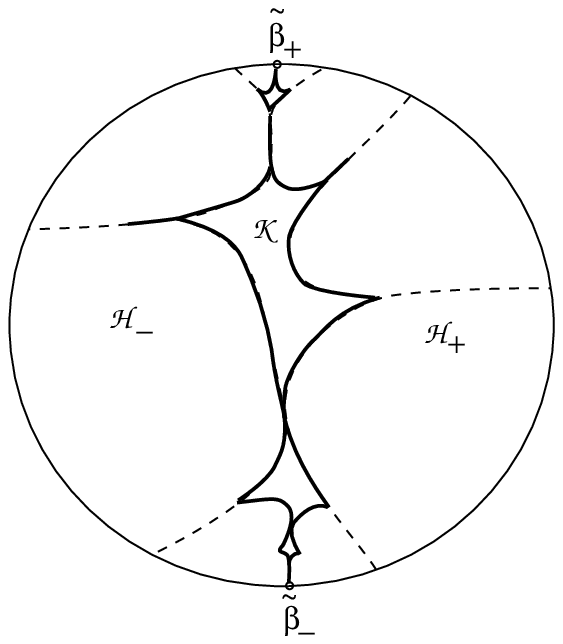}{Part of the region $\KK$, outlined
in solid lines. The dotted lines indicate a few of the halfplanes
comprising $\HH_+$ and $\HH_-$.}

It remains just to check that the projected extension is still a
train
track, i.e. that none of the diagonals cross each other. But this
is
the same as checking that, if we translate the construction
upstairs
by $\pi_1(S)$, no region is traversed by two crossing diagonals.
If this were to happen we clearly would obtain two translates of
$\til\beta$ whose endpoints separate each other, which would
contradict the assumption that $\beta$ is simple. 

This concludes the proof in the case that $\tau' = \tau$. Now in
general, note that we have so far proved the following: If
$\beta$ is not
carried in any extension of $\tau$, then there is an edge $e$ of
$\tau$ whose lift $\til e$  has the property that all train
routes of
$\til\tau$ through $\til e$ separate $\boundary\tbeta$. Now
consider a
train route $r'$ of the lift $\til\tau'$ of $\tau'$ which passes
through $\til e$. Since $\til\tau'$ is a diagonal extension of
$\til
\tau$, $r'$ must be sandwiched in between two routes of
$\til\tau$
that pass through $\til e$. It follows that $r'$ also separates
$\boundary\tbeta$. Hence when $\alpha$ is carried on $\tau'$, it
lifts
to $\alpha(e)$ routes of $\til\tau'$ through $\til e$, and 
as before we obtain $i(\alpha,\beta) \ge \alpha(e)$.
\end{pf}

\subsection{Infinite diameter and action of pseudo-Anosovs}
\label{inf diam}
We now have sufficient tools to prove the following result on the
action of $\Mod(S)$ on $\CC(S)$: 

\begin{proposition}{positive translation distance}
For a non-sporadic surface $S$ there exists $c>0$ such that, 
for any pseudo-Anosov $h\in \Mod(S)$, 
any $\gamma\in \CC_0$ and any $n\in\Z$,
$$d_\CC(h^n(\gamma),\gamma)\geq c|n|.$$
\end{proposition}

As an immediate corollary we have 
$$
\diam(\CC(S)) = \infty,$$
which gives part of the conclusion of
Theorem \ref{Hyperbolicity}.  (In fact F. Luo has pointed out an
easier proof that the diameter is infinite, which we sketch here: Let
$\mu$ be a maximal geodesic lamination and let $\gamma_i$ be any
sequence of closed geodesics converging geometrically to $\mu$. Then
if $d_\CC(\gamma_0,\gamma_n)$ remains bounded, after restricting to a
subsequence we may assume $d_\CC(\gamma_0,\gamma_n)=N$ for all $n>0$. For
each $\gamma_n$ we may then find $\beta_n$ with
$d(\beta_n,\gamma_n)=1$ and
$d(\gamma_0,\beta_n) = N-1$. But $\gamma_n\to\mu$ and $\mu$ maximal
implies $\beta_n\to\mu$ as well, since $\gamma_n$ and $\beta_n$ are
disjoint in $S$. Proceeding inductively we arrive at
the case $N=1$, and the conclusion is that $\beta_n\to \mu$ and
$\beta_n=\gamma_0$, a contradiction.)

Proposition \ref{positive translation distance} should be compared to
a property of the action of a word-hyperbolic group $G$ on its Cayley
graph $\Gamma$: namely, for a fixed $c(G)>0$, if $h\in G$ has infinite
order then 
$d(h^n\gamma, \gamma) \ge c|n|$ for all $n\in\Z$,
$\gamma\in\Gamma$ (see
\cite{gromov:hypgroups,ghys-harpe,short:notes}). 

In Gromov's terminology \cite{gromov:hypgroups}, Proposition
\ref{positive translation distance} says that the action of a
pseudo-Anosov $h$ on $\CC(S)$ is {\em hyperbolic}. In general there
are two more types of isometries of a $\delta$-hyperbolic space: {\em
  elliptic}, for which orbits are 
bounded, and {\em parabolic}, for which orbits are unbounded but
$\inf d(\gamma,h^n(\gamma))/|n| =  0$. 
When $h$ is not
pseudo-Anosov, it must be reducible or finite-order. In either case
some vertex of $\CC(S)$ is fixed by a finite power of $h$, and hence
$h$ is elliptic.
Thus it follows from Proposition \ref{positive
translation distance} that the
action of $\Mod(S)$ on $\CC(S)$ has no parabolics.

\begin{pf*}{Proof of Proposition \ref{positive translation distance}}
A pseudo-Anosov map $h:S\to S$ determines measured laminations 
$\mu,\nu\in\ML(S)$, called stable and unstable laminations,
with the following properties
(for more about pseudo-Anosov homeomorphisms, see e.g.
\cite{travaux,wpt:surfaces,bers:pseudoanosov}).
They are transverse to each other,
and the complementary regions of each are
ideal polygons or once punctured ideal polygons. 
Both projective classes $[\mu]$ and $[\nu]$ in
$\PML(S) = (\ML(S)\setminus\{0\})/\R_+$ are fixed points for $h$, 
such that $[\mu]$ is attracting in $\PML(S)\setminus[\nu]$, and
$[\nu]$ is repelling in $\PML(S)\setminus[\mu]$.
In particular, since $\mu$ and $\nu$ cannot have closed-curve
components, every vertex of $\CC(S)$ approaches $[\mu]$ under
forward iteration of $h$, and $[\nu]$ under backward iteration.

Let $\tau_0$ be a generic train-track formed from a regular
$\epsilon$ neighborhood of $\mu$.  If $\ep$ is sufficiently small, the
complementary domains of $\tau_0$  are in one-to-one correspondence with
those of $\mu$, and $\tau_0$ is birecurrent (see Penner-Harer
\cite{penner-harer}). 

One can homotope $h$ to a standard form in which it permutes the
complementary regions of $\mu$, is expanding on the leaves of
$\mu$, and contracting in the transverse direction near
$\mu$. Thus, the image train-track $h(\tau_0)$ must be carried in
$\tau_0$, and fills it.

If $\tau\in E(\tau_0)$ is a diagonal extension, then $h(\tau)$ is
carried in some $\tau'\in E(\tau_0)$ by Lemma \ref{easy nesting}.
Since the number of tracks in $E(\tau_0)$ is bounded in terms of the
topology of $S$, there is some $k_0(S)$ such that, for some $k\le k_0$
the power $h'=h^k$ takes $\tau$ to a track carried by $\tau$.
Let $\BB$ denote the branch set of $\tau$, and $\BB_0\subset \BB$ the
branch set of $\tau_0$.  In the coordinates of $\R^\BB$ we may
represent $h'$ as an integer matrix $M$, with a submatrix $M_0$ giving
the restriction to $\R^{\BB_0}$ (see Penner \cite{penner:dilatation}).
Penner shows in \cite{penner:dilatation} that $M_0^n$ has all positive
entries where $n$ is the dimension $|\BB_0|$, and in fact $|M_0^n(x_0)|
\ge 2|x_0|$ for any vector $x_0$ representing a measure on
$\tau_0$. Indeed $M_0$ has a unique eigenspace in 
the positive cone of $\R^{\BB_0}$, which corresponds to $[\mu]$.
On the other hand, for a diagonal branch
$b\in\BB\setminus\BB_0$ we have, by Lemma \ref{diagonals nesting}, that
$|M^i(x)(b)| \le m_0|x|$ for all $x\in\R^\BB$ and all powers $i>0$.
Now, any transverse measure $x$ on $\tau$ must put some positive
measure on a branch of $\BB_0$, since $\tau_0$ is generic. It follows
immediately that given any 
$\delta>0$ there exists $m_1$, depending only on $\delta$ and $S$, such
that for some $m\le m_1$ we have $\max_{b\in{\BB\setminus\BB_0}}
h^m(x)(b) \le \delta \min_{b\in\BB_0} h^m(x)(b)$, for any $x\in
P(\tau)$. Applying this to each $\tau\in E(\tau_0)$, and invoking Lemma
\ref{interior of extension}, 
we conclude that, for suitable choice of $\delta$, 
$$h^m(PE(\tau_0))\subset int(PE(\tau_0)).$$

Thus letting $\tau_j = h^{mj}(\tau_0)$ we
find by induction that $PE(\tau_{j+1})\subset int(PE(\tau_j))$.
Now if $\beta\in\CC_0(S)$ satisfies
$\beta\notin PE(\tau_0)$ but
$h^m(\beta)\in PE(\tau_0)$, we have that  $h^{km}(\beta)\in
PE(\tau_{k-1})$ for $k \ge 1$.  But then 
Lemma \ref{1-nbd is closure} applied inductively  shows that
$d_{\CC}(h^{km}(\beta),\beta) \ge k$. 

For arbitrary $n\in\Z$ we note, since $h$ is an isometry on $\CC$, that
$|n| \le d_\CC(h^{nm}(\beta),\beta) \le md_\CC(h^n(\beta),\beta)$, and
conclude that $d_\CC(h^n(\beta),\beta) \ge |n|/m$.

Finally, for arbitrary $\gamma\in\CC_0(S)$ we note that 
$[h^n(\gamma)] \to [\mu]$ as $n\to\infty$, and (see the discussion before
Lemma \ref{interior of extension}) some neighborhood of $\mu$
is contained 
in $int(PE(\tau_0))$. Thus eventually $h^n(\gamma)\in int(PE(\tau_0))$.
On the other hand $[h^{-n}(\gamma)]\to
[\nu]$ as $n\to\infty$, and $\nu$ is not contained in $PE(\tau_0)$
since by the previous discussion all of $[PE(\tau_0)]$ converges to
$[\mu]$ under iterations of $h$. Since $[PE(\tau_0)]$ is closed 
it misses a neighborhood of $[\nu]$, and we
conclude that there is 
some  $p\in\Z$ for which $h^p(\gamma)\notin
PE(\tau_0)$ but $h^{m+p}(\gamma)\in PE(\tau_0)$.  Applying the
previous two paragraphs to  $\beta = h^p(\gamma)$, we obtain the
desired bound for $\gamma$ as well.

We thus have our conclusion for $c = 1/m_1$, which by construction is
independent of $h$ or $\gamma$.
\end{pf*}

\subsection{The nesting lemma}
\label{nesting proof}
We will need the following notation. If $\alpha\in \CC_0(S)$ and 
$\sigma,\tau$ are train tracks, let $d_\CC(\alpha,\sigma)$ denote
$\min_v d_\CC(\alpha,v)$ and $d_\CC(\sigma,\tau) =
\min_{v,w}d_\CC(v,w)$, where $v$ ranges over the vertices of
$\sigma$
and $w$ ranges over the vertices of $\tau$.

The goal of this section is the following lemma, whose proof will
appear at the end of it.
\begin{lemma+}{Nesting Lemma}
There exists a $D_2>0$ such that, whenever 
$\omega$ and $\tau$ are large recurrent generic tracks and
$\omega\carriedby\tau$, if
$d_\CC(\omega,\tau)
\ge
D_2$, we have
$$
PN(\omega) \subset int(PN(\tau)).
$$
\end{lemma+}

Given a train-track $\tau$ and a measure $\mu\in P(\tau)$ we can
define a {\em combinatorial length} $\ell_\tau(\mu)$ as
$\sum_b \mu(b)$, where the sum is over the branches $b$ of
$\tau$.
Similarly if $\mu \in PN(\tau)$ we can define
$\ell_{N(\tau)}(\mu)$ as the minimum of 
combinatorial lengths in the tracks of $N(\tau)$ that carry
$\mu$.

There are some easy consequences of there being only finitely
many combinatorial types of train-tracks on $S$.
For example, if $\lambda\in P(\tau)$ and one writes $\lambda$ as
a
combination $\sum_i a_i\alpha_i$ (not necessarily unique) of the
vertices $\alpha_i$ of $\tau$ with nonnegative coefficients,
then
\begin{equation}
\label{vertex coefficients and length}
\max_i a_i \le \ell_\tau(\lambda) \le C_1 \max_i a_i
\end{equation}
where $C_1$ depends on a bound on the number of vertices of
$\tau$,
and a bound $C_0$ for $\ell_\kappa(\omega)$ over all train tracks
$\kappa$ and vertices $\omega$.

Another consequence of finiteness  is that there is a constant
$B$
depending only on $S$, 
such that any two vertices of a train-track are $\CC$-distance at
most $B$ apart. (We conjecture that $B=2$).

Furthermore we have: 
\begin{lemma}{distance implies weight}
Given $L>0$ there exists $D_0(L)$ so that, if $\alpha\in P(\tau)$
and
$\ell_{\tau}(\alpha) <L$ then $d_\CC(\alpha,\tau)<D_0$.
\end{lemma}

\begin{pf}
Fixing $L$ and $\tau$, only finitely many curves $\alpha$ are
carried
by $\tau$ with $\ell_\tau(\alpha)\le L$. Thus there is an upper
bound
on their distance from the vertices of $\tau$. Taking a maximum
over
all combinatorial types of train-tracks in $S$, we have the
desired
statement.
\end{pf}

We also observe:
\begin{lemma}{distance implies large}
If $\alpha\in P(\tau)$ and 
$d_\CC(\alpha,\tau) \ge 3$ then $\alpha$ fills a large subtrack
of $\tau$.
\end{lemma}

\begin{pf}
Suppose that $\alpha$ is carried in $\kappa\subtrack\tau$ which
is not
large. Then $S\setminus\kappa$ contains a nontrivial,
nonperipheral
curve $\beta$, so that $d_\CC(\beta,\alpha) \le 1$ and
$d_\CC(\beta,v)
\le 1$ for any vertex $v$ of $\kappa$. By the triangle inequality
$d_\CC(\alpha,v)\le 2$, and since $v$ is also a vertex of $\tau$,
$d_\CC(\alpha,\tau) \le 2$.
\end{pf}

The next lemma addresses the following issue. A closed curve
carried on an
extension of a track $\sigma$ does not necessarily trace through
any
complete cycle on $\sigma$. However, if $\sigma$ is sufficiently
deeply nested in $\tau$, then any curve on an extension of
$\sigma$
is forced to run through a cycle of $\tau$, and in fact must put
a
definite amount of weight on that cycle.

\begin{lemma}{extension hits cycle}
There exists $M_0$, and for any $L$ there exists $D_1(L)$
such that if $\sigma$ is large,
$\sigma\carriedby\tau$, and $d(\sigma,\tau) \ge D_1(L)$ then the
following holds. 
Suppose $\sigma'\in E(\sigma)$ and $\tau'\in E(\tau)$, and
$\sigma'\carriedby\tau'$. Then any curve $\beta$ carried on
$\sigma'$ can be
expressed in $P(\tau')$ as $\beta_\tau + \beta_\tau'$, where
$\beta_\tau\in
P(\tau)$, and 
\begin{equation}
\label{diagonals upper bound}
\ell_{\tau'}(\beta_\tau') \le M_0\ell_{\sigma'}(\beta),
\end{equation}
\begin{equation}
\label{main track lower bound}
\ell_\tau(\beta_\tau) \ge L\ell_{\sigma'}(\beta).
\end{equation}
\end{lemma}

\begin{pf}
It suffices to prove the lemma when $\beta$ is a vertex $v$ of
$\sigma'$. For the general case, express $\beta$ as a combination
of
vertices and use (\ref{vertex coefficients and length}).

Let $W_0$ be a bound (by finiteness) on the weights
that
$v$ puts on any branch of $\sigma'$, so that by Lemma
\ref{diagonals
nesting} $v$ puts at most $m_0W_0$ on the branches of
$\tau'\setminus\tau$.
Write the vertices of $\tau'$ as
$\{\alpha_i\}\union\{\gamma_j\}$,
where $\alpha_i$ are the ones supported in $\tau$.
Then in the coordinates of $P(\tau')$ we may write $v = v_\tau +
v_\tau'$ where $v_\tau = \sum a_i\alpha_i$ and $v_\tau' = \sum_j
c_j\gamma_j$, with $a_i,c_j\ge 0$. 

For each branch $b$ of $\tau'\setminus\tau$ we have 
$v(b) = \sum c_j \gamma_j(b) \le m_0W_0$. Since for each $j$ some
$b$ has $\gamma_j(b) \ge 1$, 
we have $c_j \le m_0W_0$. 

We have shown
$$
\ell_{\tau'}(v_\tau') \le m_0W_0C_0,
$$ where recall $C_0$ is a bound for the combinatorial length of
all vertices of a train track.  
Letting $M_0 = m_0W_0C_0$, and noting that $\ell_{\sigma'}(v) \ge
1$,
we have the first desired inequality (\ref{diagonals upper
bound}).  
Let $D_0$ be the distance bound provided by Lemma \ref{distance
implies weight} for a length bound of
$C_0L+M$.
Then let $D_1 = 2B + D_0$.

Since the distance between a vertex of $\sigma'$ and any
vertex of
$\sigma$ is at most $B$, and the same for $\tau'$ and $\tau$, 
we conclude from  the assumption
$d_\CC(\sigma,\tau) \ge D_1$ that
we have 
$d_\CC(v,\tau') \ge D_0$,
and by Lemma \ref{distance implies weight}, we have 
$\ell_{\tau'}(v) \ge C_0L+M_0$.  
Now $\ell_\tau(v_\tau) =
\ell_{\tau'}(v) - \ell_{\tau'}(v_\tau') \ge C_0L$, and since
$\ell_{\sigma'}(v) \le C_0$ we have the second inequality 
(\ref{main track lower bound}).
\end{pf}

\subsection*{Proof of Lemma \lref{Nesting Lemma}}
Let $\omega\carriedby\tau$ with $d_\CC(\omega,\tau) \ge D_2$,
where
$D_2$ will be determined shortly. Let $\sigma$ be any large
subtrack of
$\omega$. We will prove that 
$PE(\sigma) \subset int(PE(\kappa))$ for some large subtrack
$\kappa$
of $\tau$. Thus by definition we will have $PN(\omega)\subset
int(PN(\tau))$, which is the desired statement.

Let $\tau=\tau_0$ and $\cdots
\tau_2\carriedby\tau_1\carriedby\tau_0$
be a sequence of tracks 
obtained from $\tau_0$ by splitting, so that
$\sigma\carriedby\tau_j$ for each $j$. 
Let $\rho$ be the first $\tau_j$ for which
$d_\CC(\tau_j,\tau)>2$.
Note that we may assume $\sigma$ fills $\rho$: in a splitting move
determined by $\sigma$, if an edge becomes empty we use a collision
move and erase the edge. Since $\sigma$ is large, $\rho$ must be as
well.

By the properties of splitting sequences (see Section
\ref{train-track intro}),
$\tau_j$ either shares a vertex with $\tau_{j-1}$ or is a
subtrack of
a track that shares a vertex with it. Thus for any vertex $v$ of
$\rho=\tau_j$,
$d_\CC(v,\tau_{j-1})\le B$,
and it follows that 
$d_\CC(v,\tau) \le 2+2B$. Therefore
$d_\CC(\sigma,\rho) \ge D_2 - 2-2B$.

Fix now $\beta$ carried by $\sigma'\in E(\sigma)$, and let us
show 
that $\beta\in int (PE(\kappa))$ for some large subtrack
$\kappa\subtrack\tau$.  The idea will be that, by Lemma
\ref{extension
hits cycle}, $\beta$ will place definite weight on some cycle of
$\rho$, and by Lemma \ref{distance implies large} this cycle will
fill
a large subtrack $\kappa_0$ of $\tau$. On the other hand $\beta$
will place
relatively little weight on any extension branches outside
$\tau$, and
we will be able to reach our conclusion for some $\kappa$
containing
$\kappa_0$. 

Since $\sigma$ fills $\rho$, by Lemma \ref{easy nesting} there is
some
$\rho'\in E(\rho)$ carrying $\beta$. 
Fix $L_1$ (to be determined shortly), and let $C_1$ be the constant in
inequality (\ref{vertex coefficients and length}).
Lemma \ref{extension hits cycle}
implies that for sufficiently large 
$D_2$ (depending on $L_1C_1$), we
can write $\beta = \beta_\rho + \beta_\rho'$
where $\ell_\rho(\beta_\rho) \ge L_1C_1\ell_{\sigma'}(\beta)$.
Inequality 
(\ref{vertex coefficients and length}) then implies that
we can write $\beta_\rho = \sum_i
a_i\alpha_i$ where  $\alpha_i$ are vertices of $\rho$, such that
$a_1 \ge L_1\ell_{\sigma'}(\beta)$. Now applying Lemma \ref{distance
  implies large},
since $d_\CC(\alpha_1,\tau) > 2$ 
there is a large subtrack $\kappa_0$ of $\tau$ such that
$\alpha_1(b) \ge 1$ for each branch $b$ of $\kappa_0$.

Therefore $\beta(b)\geq L_1\ell_{\sigma'}(\beta)$ for every 
branch $b$ of $\kappa_0$. However, 
we don't know if $\beta\in PE(\kappa_0)$. The trouble is that the
enlargement of $\kappa_0$ that supports $\beta$ may not be a
diagonal extension. 
Thus we will find an intermediate track between $\tau$ and  
$\kappa_0$ by adding branches to $\kappa_0$ that have too much
weight to  
be pushed to the corner, and show that this process terminates
with
the desired track.

Let $\tau'\in E(\tau)$ be a track carrying $\beta$ (by Lemma
\ref{easy
nesting}). By Lemma 
\ref{diagonals nesting} we know that
$\beta(b) \le m_0\ell_{\sigma'}(\beta)$ for any branch $b$ of
$\tau'\setminus\tau$. 

If for all branches $c$ of $\tau'\setminus\kappa_0$ which meet
$\kappa_0$ we have $\beta(c) < \delta L_1\ell_{\sigma'}(\beta)$
then
by Lemma \ref{interior of
extension}, $\beta\in int(PE(\kappa_0))$ and we are done.
If not, let $c$ violate this inequality, and define an extension
$\kappa_1$ of $\kappa_0$ containing $c$ as follows: let $c_\pm$
be
the
ends of $c$ where $c_-\in\kappa_0$. If $c_+\in\kappa_0$ then
$\kappa_1 = \kappa_0\union c$ is a train track. If not, then
$c_+$
is
incoming to some switch with at most $m_1$ branches outgoing 
($m_1 = m_1(S)$). At least one of those, $c_1$, a branch of
$\tau'$, has measure $\beta(c_1)
\ge {1\over m_1}
\beta(c) \ge {\delta\over m_1} L_1\ell_{\sigma'}(\beta)$. Add
this
branch, and continue adding branches
of $\tau'$ until we find one which touches $\kappa_0$ again. Let
$\kappa_1$
denote 
$\kappa_0$ together with this chain of branches, and note that
for
all
branches $b$ of $\kappa_1$, 
$\beta(b) \ge {\delta\over m_2}L_1\ell_{\sigma'}(\beta) $,
where $m_2 = m_2(S)$. 

If now there is no edge of $\tau'\setminus\kappa_1$ adjacent to
$\kappa_1$ with measure 
at least $\delta {\delta\over m_2}L_1\ell_{\sigma'}(\beta)$, we
are
done. Otherwise, 
we can repeat this process, obtaining a sequence $\kappa_i$ of
extensions, each of which is a subtrack of $\tau'$ and  which
must terminate after at most $m_3$ steps, $m_3=m_3(S)$. 
Thus, $\beta(b)$ for any branch $b$ of $\kappa_i$ is always at
least
$\left({\delta\over m_2}\right)^{m_3} L_1\ell_{\sigma'}(\beta)$.
If we
have chosen $L_1$
sufficiently large that
$m_0 <  \left({\delta\over m_2}\right)^{m_3} L_1$,
this process must terminate {\em without} appending to $\kappa_j$
any branches of $\tau'\setminus\tau$, since as above all such branches
have $\beta$-measure at most $m_0\ell_{\sigma'}(\beta)$.
Therefore we must end with some
$\kappa_j\subtrack \tau$ for which $\beta\in int( PE(\kappa_j))$,
and we are done.
\qed
\medskip

We note that 
a corollary of the proof is the following quantitative version of
the Nesting Lemma, obtained by taking the constant $L_1$
sufficiently large:

\begin{lemma}{quantify nesting}
The constant $D_2$ in the Nesting Lemma may be chosen so that, if
$\omega\carriedby\tau$ and $d_\CC(\omega,\tau) \ge D_2$, then for
any
$\beta\in PN(\omega)$ there is a subtrack $\kappa\subtrack\tau$
such
that $\beta\in PE(\kappa)$ and, for any branch $b$ of $\kappa$,
$$
\beta(b) \ge 2 \ell_{N(\omega)}(\beta).
$$
\end{lemma}

\subsection{Growth of intersection numbers}
Lemma \ref{quantify nesting} implies, in particular, that the
combinatorial length of a curve carried on a train track grows at least
exponentially with its distance from a fixed point, say a vertex
of the track. As a consequence (see also Lemma \ref{quantitative
intersection}), its intersection number with any fixed curve not
carried on the track should grow exponentially.

A finer analysis shows that, if two curves are both far from a
fixed one and relatively near each other, then both are deeply
nested in diagonal extensions of the same track, and as a
consequence their intersection numbers with any fixed curve are
very large compared to their intersection number with each other.
In the closing argument of the Projection Theorem \ref{Projection
Theorem}, in Section \ref{projection proof}, we will use the
following quantitative version of this observation.


\begin{lemma}{nesting consequence}
Given $Q,k>0$ there exist $D_3,\nu$ such that the following
holds. 
If $\alpha,\beta$ and $\gamma$ in $\CC_0(S)$ are such that
$$ d_\CC(\beta,\alpha) \ge D_3$$
and
$$ d_\CC(\gamma,\beta) \le \nu d_\CC(\beta,\alpha) $$
then
$$
\min_{\alpha'} i(\beta,\alpha') \cdot
\min_{\alpha'} i(\gamma,\alpha')
\ge
Q i(\beta,\gamma),
$$
where $\alpha'$ varies over the 
$k$-neighborhood of $\alpha$ in $\CC$.
\end{lemma}

\begin{pf}

Extend $\alpha$ to a pair of pants decomposition of $S$.  Such a
decomposition determines a family of {\em standard train-tracks},
obtained by choosing one of a finite number of configurations in
each pair of pants and in a connecting collar between any
adjacent pairs of pants.  (See Penner-Harer \cite{penner-harer}).
Each such track is generic and birecurrent, and
The family of tracks has the property that {\em any} simple
closed curve on $S$ is carried by one of them. 
Thus, let
$\tau_0$ denote a standard train-track carrying $\beta$.
Depending on the choice of local picture in an annulus
neighborhood of $\alpha$, 
\marginpar{figure?}
$\alpha$ is either a vertex cycle of $\tau_0$, or is distance at
most 2 from a vertex cycle. 
It follows that $d_\CC(\beta,\tau_0)
\ge d_\CC(\beta,\alpha) - (2+B)$.

Since $\beta$ is assumed far from $\alpha$, by Lemma
\ref{distance implies large} it
fills a large subtrack of $\tau_0$, which we will continue to
call $\tau_0$.

Now we will show that we can find a sequence
$\tau_n\carriedby\cdots\carriedby\tau_0$ of train-tracks,
each carrying $\beta$, so that $d_\CC(\tau_{j+1},\tau_j) > D_2$
and the length of the sequence is $n \ge d_\CC(\beta,\tau)/(D_2 +
2B)$.  This is done by splitting. Perform a sequence of
splittings of $\tau_0$ determined by the weights of $\beta$, and
terminating with a track 
that has $\beta$ as a vertex.
Inductively define $\tau_{j+1}$ to be the first track in the sequence
such that $d_\CC(\tau_{j+1},\tau_j) > D_2$. 
Since $\tau_{j+1}$ is
the first,  we may conclude
as in the proof of Lemma \ref{Nesting Lemma} that 
$d_\CC(v,\tau_j)\leq D_2+2B$ for any vertex $v$ of $\tau_{j+1}$.
It follows that $d_\CC(\beta,\tau_{j+1}) \ge
d(\beta,\tau_j) - (D_2 + 2B)$, and we may continue.

Lemma \ref{Nesting Lemma} now guarantees that $PN(\tau_{j+1})
\subset
int(PN(\tau_j))$. 

If $d_\CC(\alpha,\alpha') \le k$ then $d_\CC(\alpha',\tau_0) \le
k+2$,
and by applying Lemma \ref{1-nbd is closure} inductively we see
that $\alpha'$ cannot be in $PN(\tau_{k+3})$. (Note that to apply
Lemma \ref{1-nbd is closure} we need each $\tau_j$ to be transversely
recurrent, but this follows from the fact that $\tau_0$ is.)

Another application of  Lemma \ref{1-nbd is closure}
shows that if $d_\CC(\beta,\gamma) \le m < n$, then
$\gamma \in PN(\tau_{n-m})$. 

Thus, assuming $n$ is sufficiently large compared to $m$, we may
conclude that both $\beta$ and $\gamma$ are contained in
$PN(\tau_{k+3})$. Furthermore,  applying Lemma 
\ref{quantify nesting} repeatedly, we also have that for a large
subtrack $\kappa$ of $\tau_{k+3}$, $\gamma$ puts weight at least 
$2^{n-m-k-3} \ell_{N(\tau_{n-m})}(\gamma)$ on every branch of
$\kappa$. The same holds for $\beta$ and a large subtrack
$\kappa'$ of
$\tau_{k+3}$. Applying Lemma \ref{quantitative intersection}, we
find that 
$$
i(\alpha',\gamma) \ge 2^{n-m-k-3}
\ell_{N(\tau_{n-m})}(\gamma),
$$
and similarly for $\beta$. 

On the other hand it is easy to see that
$$
i(\beta,\gamma) \le C_2
\ell_{N(\tau_{n-m})}(\beta)\ell_{N(\tau_{n-m})}(\gamma). 
$$
where $C_2$ depends only on the topological type of $S$. This is
because in every branch of $\tau_{n-m}$ a strand of $\beta$ and
one of
$\gamma$ can only have one essential intersection, and a strand
in a diagonal
branch of $N(\tau_{n-m})$ can only hit strands in diagonal
branches, 
and two diagonal branches can intersect at most once if the
complementary domain is a polygon, or twice if it is a punctured
polygon. 

Putting these inequalities together, if $n-m-k$ is sufficiently
high we have the desired inequality.
\end{pf}

\section{Geometry of quadratic differentials}
\label{qd geom}

Let $q$ be a holomorphic quadratic differential of area 1 with
respect to some conformal structure $x$ on $S$. In this section
we will study the geometry imposed by $q$, with particular regard
to the way nearly horizontal and nearly vertical geodesics are
arranged, and how they intersect each other.  Our main goals are
the Vertical Domain Lemma
\ref{Vertical domain properties}, which gives a particular
``thickening'' of a nearly vertical curve with some useful properties,
and the Intersection Number Lemma \ref{beta and gamma intersect},
which gives conditions for a nearly vertical and a nearly horizontal
curve to have large intersection number.

\subsection{Basic properties and uniform estimates}
\label{basic qd}
A {\em straight segment} with respect to $q$ is a path containing
no singularities in its interior, and which is geodesic in the
locally Euclidean metric of $q$.  If $S$ has no punctures, 
a {\em geodesic segment} is
composed of straight segments which meet at singularities making
an angle of at least $\pi$ on either side. A straight segment
connecting two singularities is also called a {\em saddle
connection}.

A {\em metric cylinder} in $q$ is an annulus which is isometric
to the product of a circle and a line segment.

When $S$ has no punctures, each nontrivial homotopy class has a
geodesic representative. However when there are punctures the
metric of $q$ is incomplete and we must slightly generalize the
notion. From now on by ``geodesic representative'' of a closed
curve $\alpha$ we mean a curve $\alpha^*$ in the compactified
surface $\hat S$ (adding the punctures) such that
$\alpha^*\intersect S$ is composed of geodesic arcs, and there is
a homotopy from $\alpha$ to $\alpha^*$ which until the last
moment is contained in $S$. One can formalize this notion by
first excising from $S$ open $r$-neighborhoods of the punctures, 
considering geodesic representatives in the resulting compact surface,
and then letting $r$ tend to $0$.
It is not hard to see that any
non-peripheral homotopy class has such a geodesic representative,
which has minimal length, and the representative is unique unless
there is a metric cylinder foliated by curves in the homotopy
class (however the homotopy class is not uniquely determined by the
representative).  The same discussion works for the geodesic
representative of a homotopy class of paths rel endpoints.

If we start with a homotopy class of simple curves then the geodesic
representative does not have to be simple: even in the absence of
punctures, it may have self-tangencies along saddle connections,
because the metric is not smooth at the zeros of $q$. Furthermore,
when a geodesic representative passes through a puncture, since the
total angle around the puncture is $\pi$, it is easy to see that in
fact the path approaches the puncture along a straight segment, and
then retraces the same segment in the opposite direction.

\subsection*{Topological constants.}
For later reference, $n_1,\ldots,n_5$ will denote the following bounds,
which may easily be computed in terms of the genus and number of
punctures of $S$. Let $n_1$ bound
the number of singularities of $q$, including
punctures.
Let $n_2$ bound the number of disjoint saddle connections
which may
appear simultaneously in $S$. 
Let $n_3$ be an isoperimetric constant, such that $\Area_q(X)
\le n_3
\diam_q(X)^2$ for any subset $X$ of $S$.
Let $n_4$ bound the size of a sequence $X_1\subset \cdots
\subset X_{n_4} \subset S $ for which $i_*\pi_1(X_j)$
is a proper subgroup 
of $i_*\pi_1(X_{j+1})$, where $i_*$ is the map induced on $\pi_1$
by inclusion
into $S$.
Let $n_5$ bound $1/\pi$ times the sum of cone angles over all
singularities of $q$. 

\subsection*{Definite collars}
Let the {\em width} of an annulus $A$ in a metric $q$ denote the
minimal distance between boundaries, and the {\em circumference}
the minimal length of a curve going once around $A$.
A compactness argument using the moduli space of Riemann surfaces
yields the following:
\begin{lemma}{Uniform collar}
There exists $W>0$ depending only on the topology of $S$ such
that, for each unit area quadratic differential $q$ there exists
a nonperipheral annulus of width
$W$.  Furthermore, given $\mu$ there exists $L>0$ so that the
annulus can be chosen either to be a metric cylinder of modulus
at least $\mu$, or to have circumference at least $L$.
\end{lemma}

\begin{pf}
If the statement is false, then there is a sequence of conformal
structures $x_i$ on $S$, unit-area holomorphic quadratic
differentials $q_i$, and $L_i\to 0$, $W_i\to 0$ such that there
is no annulus in $(S,q_i)$ of width at least $W_i$ which either
has circumference at least $L_i$ or is a metric cylinder of
modulus $\mu$.

We can now apply a compactification argument whose details may be
found in Masur \cite{masur:iet}.  We may take a subsequence so
that $(S,x_i)$ converge in a compactified moduli space to a noded
Riemann surface $(S',x)$, where $S'$ may be taken as the
complement in $S$ of a collection of disjoint curves, and $q_i$
converge on compact sets of $S'$ to some $q$. Given $\mu>0$ there
is a $K(\mu)>0$ (depending on the topological type of $S$) such
that the following alternative holds: if $\diam(q_i) \ge K$ for
all sufficiently high $i$ then eventually $(S,x_i,q_i)$ contains
a metric cylinder of width at least $W$, and modulus at least
$\mu$.  In this case we have contradicted the choice of sequence,
hence we are done.  If $\diam(q) \le K$ for all sufficiently high
$i$ then the limiting $q$ is non-zero on at least one component
$R$ of $S'$. (The two possibilities are not mutually exclusive).
We also note that $q$ has at most simple pole singularities at
the punctures.

Since $R$ supports a non-zero holomorphic quadratic differential
of finite area, it cannot be a sphere with less than 4 punctures.
It follows that there is some simple, nontrivial, nonperipheral
curve in $R$, so let $A$ be any collar for this curve. 
If $W$ and $L$ are the width and circumference of $A$, then in
the approximating metrics of $q_i$ for high enough $i$ we obtain
annuli  of nearly these width and circumferences, again
contradicting the choice of sequence.
\end{pf}

\subsection*{Definite boxes}
We will need the following notion, where a {\em rectangle}
denotes an embedded
Euclidean rectangle with respect to $q$, in particular containing
no
singularities in its interior.
\begin{definition}{Box definition} 
Let $\omega$ denote a $q$-geodesic segment or closed curve.
If $N,\delta>0$, an $(N,\delta)$ box for $\omega$
is a rectangle containing at least $N$
parallel strands of $\omega$ (counting multiplicity)
of equal length $\delta$, parallel to two
of the sides of the rectangle.  The endpoints of the strands are
on the orthogonal sides of the rectangle. The lengths of the
orthogonal sides are at most $\delta$.
\end{definition}
{\bf Remark:} Note that if $N<1$, an $(N,\delta)$ box means a 
$(1,\delta)$-box.

As a consequence of Lemma \ref{Uniform collar} we can prove the
following:
\begin{lemma}{Box}
Let $q$ be a unit-area holomorphic quadratic differential on
$(S,x)$, and suppose that there are no $q$-metric cylinders of
modulus greater than 2 in $S$. Let $A$ denote the nonperipheral
annulus of width $W$ and length $L$ provided by Lemma
\ref{Uniform collar}.  There exist $\delta,r>0$, depending only
on the topology of $S$, such that for any closed $q$-geodesic
$\gamma$ which has intersection number $N>0$ with the core of
$A$, there is a $(rN,\delta)$-box for $\gamma$ in
$A$. Furthermore, the $q$-injectivity radius at the center of the
box is at least $\delta$.
\end{lemma}

\begin{pf}
On a smaller annulus $A'\subset A$ of width $W/2$ the
$q$-injectivity radius is at least $\delta_1 =
\min(W/4,L/2)$. There are $N$ segments (with multiplicity) of
$\gamma$ of length $W/2$ passing through $A'$.  Centered on any nonsingular
point of a segment $\sigma$ of $\gamma \intersect A'$ there is a
geodesic segment orthogonal to $\sigma$ of length $2\delta_1$,
and so (recalling $n_1$ from above) there must be a segment on
$\sigma$ of length at least 
$W/2n_1$ for which these orthogonal segments meet no
singularities, and therefore make a $(1,W/2n_1)$ box for $\sigma$,
with center on $\sigma$

Consider all such boxes in $A$. There are $N$ (with multiplicity)
and we must check that there is sufficient overlap. For each box
consider a box of half the size with the same center.  Since each
box has definite area and the area of $q$ is 1, we find that
there must be a point simultaneously in $\max(rN,1)$ half-boxes
for a fixed $r>0$, and hence a box containing $\max(rN,1)$
centers of boxes. This is the desired box.
\end{pf}

\subsection{Vertical and horizontal}
\label{vertical and horizontal} From now on, let us suppose that two
constants $\theta,\ep>0$ have been fixed 
satisfying a short list of constraints which will appear in the
course of the proof. For now assume $\theta<\min(1/2,\ep^2)$.

\begin{definition}{almost vertical def}
We say a straight segment is {\em almost vertical} (respectively
{\em almost horizontal}) with respect to $q$ if its direction is
within $\theta$ of the vertical (resp. horizontal) direction of
$q$.  We say a geodesic segment or closed curve is almost
vertical (resp.  almost horizontal) if it is composed of straight
segments each of which is almost vertical (resp. almost
horizontal) or has length at most $\epsilon$.
\end{definition}
Note that a (weak) consequence of the condition $\theta<1/2$ is
that 
$$ |\alpha|_{q,v} > \half |\alpha|_q.$$

We now define a certain type of thickening, which we call a {\em
vertical} (or {\em horizontal}) {\em domain}, that will be useful
in several places.  
\begin{definition}{vertical domain def} 
Let $\omega$ be an almost vertical geodesic segment or closed
curve.  The {\em vertical domain} $\Omega_\ep(\omega)$ is
constructed as follows.  For any  point $p\in\omega$ let
$\sigma_p$ be the maximal open $q$-horizontal segment about $p$ which
contains no singularities or punctures, and such that each
component of $\sigma_p-\{p\}$ has length at most $\ep$. Let
$\Omega_\ep(\omega)$ be the closure of $\union_{p\in\omega}
\sigma_p$.

We similarly define a {\em horizontal domain} $\Psi_\ep(\omega)$
if $\omega$ is almost horizontal, where the $\sigma_p$ are
vertical segments.
\end{definition}

\realfig{vertical domain sketch}{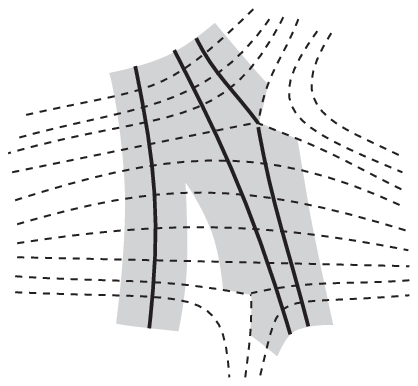}{An example of a
vertical domain. The horizontal foliation is dotted, $\omega$ is solid
and $\Omega_\ep(\omega)$ is in grey}


Let us record some useful properties of this construction.

\begin{lemma}{Vertical domain properties}
There exist positive $L_0$  and $a_0,a_1,a_2$
such that the following holds. Let $q$ be a unit-area quadratic
differential on $S$ and let $\omega$ be an almost vertical
geodesic representative of a 
simple closed curve or segment (rel endpoints).
If $\omega$ is a segment, assume it has length
at least $L_0$.
Let $\tau$ be an almost-horizontal straight segment of diameter 
$\diam_q(\tau)\ge a_0 \ep$, 
which is disjoint from $\omega$. Let $\Omega =
\Omega_\ep(\omega)$.
Then we have:  
\begin{enumerate}
\item \label{non-simple}
The map $\pi_1(\Omega)\to\pi_1(S)$ induced by inclusion has 
non-trivial, non-peripheral image.
\item \label{leave out tau}
There is a subsegment of $\tau$ of diameter at least
$a_1\diam_q(\tau)$
which is disjoint from $\Omega$.
\item  \label{boundary length}
The boundary of $\Omega$ has $q$-length
at most $a_2\ep + 1/\ep$
\item \label{horizontal boundary length}
The boundary of $\Omega$ has 
horizontal length at most $a_2\ep$.  
\end{enumerate}
\end{lemma}

\begin{pf} 
Let $n_1$,  $n_2$ and $n_5$ be the topological constants described in
\S\ref{basic qd}, and choose $L_0  = 2n_2\ep + 6/\ep$. We will see that
in fact the $a_i$ can be written explicitly in terms of the $n_i$.

We first prove Part (\ref{non-simple}).
If $\omega$ passes through a puncture then 
it follows from the definition that
$\Omega$ contains a neighborhood of the puncture. Thus if $\omega$ is
the (generalized) geodesic representative of a closed curve, or even
if it passes through punctures more than once, 
part (\ref{non-simple}) is obviously satisfied. Therefore, possibly
restricting to a subsegment, we may assume that $\omega$ is a segment
of length at least $L_0/2$, passing through no punctures in its
interior. 

For all non-singular
$p\in\omega$ let $\sigma_p$ be the horizontal arcs of Definition
\ref{vertical domain def}. For all $y\in \Omega$ let 
$f(y) = \#\{p:y\in\sigma_p\}$ (where the number is counted with
multiplicity if $\omega$ is not embedded).
Then $\Omega$ can be described as the closure of a finite
union of open parallelograms with horizontal sides of
length $2\ep$, whose heights sum
to $|\omega|_{q,v}$, and $f$ gives the degree of overlap of
these parallelograms.  It follows that $\int_\Omega
f(y) = 2\ep|\omega|_{q,v}$ (where the integral is with respect to
$q$-area). On the other hand 
$\int_\Omega f(y) \le \max(f)\Area(\Omega)$.
Since the area of $q$ is 1 we conclude 
$$
        \max f \ge {2\ep |\omega|_{q,v}}.
$$
If $\omega'$ is the union of almost vertical straight segments of
$\omega$ and $\omega''$ is 
the rest, then $|\omega''|_q \le n_2\ep$, so $|\omega'|_q \ge L_0 -
n_2\ep$. Since $|\omega'|_{q,v} \ge \half|\omega'|_q$, 
the choice of $L_0$ guarantees that $\max f \ge 3$.

We conclude that there is a point $y$ contained in a
horizontal segment $\sigma\subset \Omega$ which cuts $\omega$ in
three places (counted with multiplicity if $\omega$ is not
embedded). If there are two consecutive 
intersection points 
on $\sigma$ with the same orientation,
then a segment  of $\omega$ together with an interval of
$\sigma$ make a simple curve in $\Omega$, geodesic except at the
intersection points  where the total turning angle (measured
in the $q$ metric) is at most $2\theta < \pi$, and hence it
cannot be trivial or peripheral. (A curve bounding a disk has
total turning angle at least $2\pi$, and a curve bounding a
puncture has total angle at least $\pi$).
If every two consecutive
intersections have opposite orientations,
we can take three consecutive points so that the orientation
matches on the outer two, produce a curve passing through a
segment of $\sigma$ and $\omega$ with one self-intersection
point, and do surgery to get a simple curve with total turning
angle at most $4\theta< \pi$. Again it must be non-trivial and
non-peripheral. This proves part (\ref{non-simple}).

We now consider parts (\ref{boundary length}) and
(\ref{horizontal boundary length}).  For any
$y\in\boundary\Omega$, there is some $z\in\omega$ joined to $y$
by a horizontal arc $\sigma$ which meets $\omega$ only in $z$.
Note that $\sigma$ may pass through a singularity or a puncture
(possibly $z$ itself). If so, then a portion of $\sigma$ may
lie on $\boundary \Omega$, and contributes at most $\ep$ to
$|\boundary\Omega|_{q,h}$. The number of horizontal arcs issuing from
a singularity is $1/\pi$ times its cone angle,
so this contribution to $|\boundary\Omega|_{q,h}$
is bounded by $n_5\ep$.

If $\sigma$ meets no singularities then $y$ is contained in a
segment of $\boundary\Omega$ parallel to a segment of $\omega$
containing $z$.  The total length of such portions of
$\boundary\Omega$ which are not almost vertical is therefore
bounded by $2n_2\ep$: there are at most $n_2$ such segments in
$\omega$, they can be approached from either side, and each has
length at most $\ep$.

Finally for any segment $\kappa$ of the portion of $\boundary
\Omega$ which is almost vertical, we note that the segments
$\sigma$ form
an embedded parallelogram of width $\ep$ and height
$|\kappa|_{q,v}$.  Since $q$ is unit area, we conclude that the
vertical length of this portion of the boundary is at most
$1/\ep$. Its horizontal length is bounded by $(1/\ep)\tan\theta< 2\theta/\ep$
(assuming $\theta<1/2$), which is bounded by $2\ep$ since
$\theta<\ep^2$.  
Putting these together we have a bound 
$|\boundary\Omega|_{q,h} \le (n_5+ 2n_2+2)\ep$, 
and 
$|\boundary\Omega|_{q} \le (n_5+ 2n_2+2)\ep + 1/\ep$, 
which proves parts (\ref{boundary length}) and
(\ref{horizontal boundary length}).

Finally it remains to prove (\ref{leave out tau}).  For $y\in
\tau\intersect \Omega$, let $\sigma_y$ be the horizontal segment
of length at most $\ep$ joining $y$ to $\omega$. Suppose that $y$
is at least $2\ep$ away from any singularity of $q$, endpoint of
$\omega$ or $\tau$, or segment of $\omega$ that is not almost
vertical.  Then (with the assumption
$\theta<1/2$) a segment of the almost-horizontal $\tau$ of length
$2\ep$ must intersect the
almost-vertical segment of $\omega$ passing through the endpoint of
$\sigma$, but we have assumed $\tau\intersect\omega = \emptyset$.
Thus, $\tau\intersect\Omega$ is contained in a
$2\epsilon$-neighborhood of the singularities, endpoints and
non-almost-vertical segments of $\omega$. There are at most $k =
n_1 + n_2 + 2$ of these, each of diameter at most $\ep$. Let $d$
be the largest diameter of a component of $\tau\setminus\Omega$.
Then the diameter of $\tau$ is bounded by $5\ep k + (k+1) d$ (the
$5\ep$ bounds the diameter of a $2\ep$ neighborhood of a segment of
diameter $\ep$).
For $\diam_q(\tau) \ge 10\ep k$, say, we find that
$d \ge \diam_q(\tau)/2(k+1)$, which gives part (\ref{leave out tau}).
\end{pf}

Let us also note the following observation which will be used in
the
proofs of Lemmas \ref{beta and gamma intersect} and \ref{Segment
Disjoint From Horizontal}.

\begin{lemma}{horizontal segment in disk}
If $q$ is a  holomorphic quadratic differential on $S$ and 
$\tau$ is an embedded straight segment in a disk $D$ in $S$,
then
$|\boundary D|_{q,h} \ge 2|\tau|_{q,h}$. If $\tau$ is in a
once-punctured disk $D'$ 
and $q$ has at most a simple pole at the puncture then
$|\boundary D'|_{q,h} \ge |\tau|_{q,h}$.
\end{lemma}
\begin{pf}
Consider first a disk $D$. For any $p\in\tau$ extend a vertical
segment $\sigma$ in both directions until it hits $\boundary D$.  This
will happen since $D$ is simply connected. It follows immediately that
the horizontal length of the portion of $\boundary D$ cut off by these
segments is equal to twice the horizontal length of $D$.  
For a punctured disk $D'$, note that a segment
$\sigma$ could hit $\tau$ at both ends, on the same side of $\tau$, if
it goes around the puncture. However since the puncture is at most a
pole this can only happen on one side of $\tau$, and the vertical
segments extended from the other side still give the desired bound.
\end{pf}

\subsection{The intersection number lemma}
\label{vertical and horizontal intersection}
Our first application of the vertical domain and box lemmas is the
following lemma, which states that an almost horizontal and an almost
vertical curve which intersect every bounded curve a definite amount
also intersect each other proportionally.
Compare this fact with Lemma \ref{nesting consequence}; the two will
be applied together to yield a contradiction.

\begin{lemma}{beta and gamma intersect}
Suppose that 
$ \ep < \min ( \delta/4,
a_1^{n_4}\delta/2a_0,a_1^{n_4}\delta/2n_4a_2)$,
in addition to previous constraints.
There exist $M,h>0$ such that, given
$x\in\TT(S)$ and unit area quadratic differential $q$ on $(S,x)$,
if $\beta$ is an almost horizontal closed geodesic 
and $\gamma$ is an almost vertical closed geodesic 
with respect to $q$, and for every
non-peripheral simple closed  $\alpha$ of $q$-length at most $M$,

\begin{equation}
\label{intersect alpha}
i(\beta,\alpha) \ge B\ \text{and} \ i(\gamma,\alpha) \ge C, 
\end{equation}
then
$$ i(\beta,\gamma) \ge hBC $$
\end{lemma}

\begin{pf}
Apply lemma \ref{Uniform collar} to get constants $L, W>0$ such that
either $q$ has a flat cylinder of modulus 2, or an annulus of
circumference at least $L$ and radius $W$.  
Let $M=\max(1/\sqrt2, 1/W,
n_4 (a_2 \ep + 1/\ep))$.
Consider now both possible cases.

\bfheading{Case A.} If $q$ has a flat cylinder of modulus 2, 
let  $\alpha$ be the core of this cylinder. 
Then $\alpha$ has $q$-length at most
$1/\sqrt{2} $, so that (\ref{intersect alpha}) gives $B$ and $C$
strands of $\beta$ and $\gamma$, respectively, crossing the
annulus.
It is easy to see that any two nearly orthogonal
segments cutting through the annulus must intersect at least
once. It
follows that $i(\beta,\gamma) \ge BC$, so with $h \le 1$, we are
done. 

\bfheading{Case B.} If $q$ has an annulus $A$ with circumference
at
least $L$ and width $W$, note that $A$ has modulus at least $W^2$
(since its area is at most 1),
and hence $Ext_x(\alpha) \le 1/W^2$, where $\alpha$ is the core
of $A$.
In particular  $|\alpha^*|_q \le 1/W$ where $\alpha^*$ is the
$q$-geodesic representative (see \S\ref{teich defs}).
Thus $\gamma$ contains
at least $C$ segments (with multiplicity) 
crossing $A$ and hence of length at least $W$,
and similarly $\beta$ contains at least $B$ such segments.

Lemma \lref{Box} guarantees an almost horizontal $(c
B,\delta)$-box
$H$ for $\beta$, with injectivity radius at least $\delta$ at its
center. Let $\tau_0$ denote the segment of length $\delta/2$
centered
on the center of $H$, parallel to the direction of $\beta$.
It has the property (recalling $\theta<1/2$) that any
almost-vertical segment that 
meets $\tau_0$ must cut through all the $\beta$ strands in $H$.

Since $\gamma$ must have length at least $CW$, we may divide it
into at least
$ CW/L_0 - 1$ pieces, each of which has
length at least $L_0$.
(If $CW/L_0 <2$ we may instead take
the whole closed curve $\gamma$, and prove the theorem for $C=1$.
The
discrepancy is absorbed in the constants.)
Each of these is almost vertical, though they may traverse saddle
connections of length at most $\ep$ which are not almost
vertical.
However, these 
short segments cannot meet $\tau_0$, since $H$ contains no
singularities (here we are using the assumption
$\ep<\delta/4$). Thus, for each segment that meets
$\tau_0$ we obtain $rB$ essential intersections with $\beta$. If
all of
them do meet $\tau_0$, then we are done. 

Thus suppose that one segment $\omega_1$ is disjoint from
$\tau_0$. 
Let $X_1 = \Omega_\ep(\omega_1)$. Lemma \ref{Vertical domain
properties} guarantees that $X_1$ generates a nontrivial,
nonperipheral subgroup of $\pi_1(S)$.
Note that the diameter of $\tau_0$
is $\delta/2$, by the injectivity
radius lower bound in $H$. Thus part (\ref{leave out tau}) of Lemma
\ref{Vertical domain properties} guarantees that a subarc
$\tau_1$ of $\tau_0$, 
with diameter $a_1\delta/2$, is disjoint from $X_1$
(to apply the Lemma we need $\delta/2 \ge a_0\ep$, which is implied by
the conditions on $\ep$).

Because the $q$-length of any nontrivial nonperipheral component of
$\boundary X_1$ is bounded by $a_2\ep + 1/\ep \le M$,
$\gamma$ intersects it essentially, $C$ times.  Thus for any component
$Y$ of $S\setminus X_1$ which is not a disk or punctured disk, there
are $C$ arcs of $\gamma$ (with multiplicity) passing through $Y$ with
both endpoints on $\boundary Y$, which are not deformable back into
$X_1$.
 

Apply this where $Y$ is the component of $S\setminus X_1$
containing
$\tau_1$. This cannot be a disk or punctured disk, because the
horizontal length of $\boundary Y$, which is at most $a_2\ep$ by
Lemma \ref{Vertical domain properties}, is
smaller than the length $a_1\delta/2$ of 
$\tau_1$ by our assumptions on $\ep$, and we may apply Lemma
\ref{horizontal segment in disk}.

Thus, if all $C$ arcs in $Y$ meet $\tau_1$ then as before we have our
required intersections between $\beta$ and $\gamma$, and we are done.

If one arc $\omega_2$ is disjoint from $\tau_1$ then define $X_2 = X_1
\union \Omega_\ep(\omega_2)$.  We may apply Lemma
\ref{Vertical domain properties} and the same arguments as before
to find a subarc $\tau_2$ of $\tau_1$, of diameter $a_1^2\delta/2$,
which is disjoint from $\Omega_\ep(\omega_2)$, and hence from $X_2$.

We may continue by induction, generating a sequence $X_1\subset\cdots
X_j\subset X_{j+1}$ and subarcs $\tau_j\subset\tau_1$ of length
$a_1^j\delta/2$ disjoint from $X_j$. At each step, $|\boundary X_j|_q$
is incremented by at most $a_2\ep+1/\ep$, and $|\boundary
X_j|_{q,h}$ goes up by at most $a_2\ep$.  Since $X_{j+1}$ cannot be
deformed into $X_j$ the process must terminate within $n_4$ steps.  By
our assumption on $\ep$ we can apply Lemma \ref{horizontal segment in
disk} each time so that the component of $S-X_j$ containing $\tau_j$
is never a disk or punctured disk.  It follows that the only way the
process can terminate is by giving $C$ intersections of $\gamma$ with
$\tau_j$ for some $j\le n_4$, which concludes the proof.
\end{pf}

\section{Proof of the projection theorem}
\label{projection proof}

In this section let  $q$
be a quadratic differential of area 1 with respect to a conformal
structure $x$ on $S$.  
let $L_q$ denote the corresponding Teichm\"uller geodesic.  
We will denote the Riemann surfaces along $L_q$ by
$x_t = L_q(t)$ where $t$ is arclength, and the 
quadratic differentials by $q_t$.

Recall that the geodesic gives rise to a map $F_q:\R\to\CC$, and
a projection $\pi=\pi_q:\CC\to \R$ defined as in section
\ref{outline}.  To prove the Projection Theorem \ref{Projection
Theorem}
we must show that this projection satisfies the contraction
property
(Definition \ref{contraction property}).

We will also assume that our constants $\ep,\theta$ satisfy the
assumptions of the Intersection Number Lemma \ref{beta and gamma
intersect}. 

\subsection{Bounded adjustments}
We will first need to examine transitions along $L_q$ from mostly
vertical to balanced to mostly horizontal curves. As measured by
the Teichm\"uller length parameter, a nearly vertical curve can
take a very long time to become balanced. However we find that in
a number of crucial situations the transition takes bounded time
(independent of $q$)
{\em as viewed in the curve complex} (that is, when considering
quantities such as $\diam_\CC(F[s,t])$ instead of $|s-t|$).

The relevant insight is illustrated by this sketch of the proof
of Lemma \ref{Segment Disjoint From Horizontal}: Consider a very
long nearly vertical segment with respect to $q_0$, which does
not fill the whole surface (say it avoids a definite-length
horizontal segment).  Then if for $t>0$ the segment is still long
and nearly vertical, it fills up some proper subsurface of $S$
which can only shrink as $t$ increases. The boundaries of the
resulting sequence of surfaces form a bounded-length sequence in
$\CC(S)$.  This is made precise using the Vertical Domain
construction.

\medskip

Our first observation about the map $F$ is that it is, on a large
scale, Lipschitz:
\begin{lemma+}{Lipschitz}  There exist $C,D>0$ such
that for any $q$, $t_1$ and $t_2$ we have
$$d_\CC(F_q(t_1),F_q(t_2))\leq C|t_2-t_1|+D.$$
\end{lemma+}
\begin{pf}  As in Lemma \ref{Phi diam bound}, let
$e_0(S)$ be such that for any
conformal structure on $S$ there is a curve with extremal length
at most $e_0$.  Suppose that $|t_2-t_1| \le 1$. Let $\alpha_i$ be
a curve of shortest extremal length for $x_{t_i}$, for
$i=1,2$. Then $Ext_{x_{t_1}}(\alpha_2) \le e^2 e_0$ (by
(\ref{kerckhoff theorem})).
A bound on
$d_\CC(\alpha_1,\alpha_2)$ follows from Lemma \ref{short are
close}.
\marginpar{Actually this is the easier side of Kerckhoff's thm,
i.e. Grotzch
inequality or something. Anyway, I put in the reference.}

The case where $|t_2-t_1|>1$ follows by subdividing.
\end{pf}

In the exceptional cases of projecting curves that are entirely
horizontal or vertical, we observe the following:
\begin{proposition}{Distance 1} If $\beta\in \CC\setminus
\CC_b(q)$
then
$d_\CC(\beta, F_q(\pi_q(\beta)))\leq 1$.
\end{proposition}
\begin{pf}
Assume that $\beta$ is vertical. Let $\Sigma$ denote the union of
compact singular leaves of the vertical foliation of  $q$, and
let
$\Sigma_\ep$ denote a regular neighborhood of $\Sigma$.  Then it
is
not hard to see (e.g. \cite{minsky:2d}) that for any
non-peripheral
curve $\gamma$ in 
$S$ the extremal length $Ext_{x_t}(\gamma)$ remains bounded as
$t\to
+\infty$ if and only if $\gamma$ can be deformed into $\Sigma$,
and 
$Ext_{x_t}(\gamma)\to 0$ as $t\to +\infty$ if and only if
$\gamma$ is
homotopic to a boundary component of $\Sigma_\ep$.

It follows that $F_q(+\infty) = F_q(\pi_q(\beta))$ is one of
these
boundary components, and since $\beta$ is in $\Sigma$, we obtain
$d_\CC(\beta, F_q(\pi_q(\beta)))\leq 1$. If $\beta$ is horizontal
we make
a similar argument, reversing the $t$-direction.
\end{pf}

In Lemmas \ref{Short Curves Become Balanced}-\ref{Almost Vertical} we
will find a series of constants $d_1$-$d_5$ 
that depend only on the topology of $S$ and not on the particular
$q$. 
The next lemma shows that the image under the map $F$ of the set
of
$t$ where a curve $\alpha$ is close to its minima has bounded
diameter
in $\CC(S)$.  In particular, once the $q_t$-length of $\alpha$ is
sufficiently short, we only need to wait a bounded amount until
it
starts to grow again.

\begin{lemma}{Short Curves Become Balanced}  There exist
$\ep_1,d_1>0$, depending only on the topology of $S$, with the
following property. If $\alpha$ is a closed $q$-geodesic
homotopic to 
a simple curve, let 
$$
     J = \{t : |\alpha|_{q_t}\le \epsilon_1\}.
$$ 
Then $\diam_\CC(F(J)) \leq d_1$.

(Note that $J$ is a bounded interval in $\R$ unless $\alpha$ is
completely vertical or horizontal.)
\end{lemma}
\begin{pf}  
Let $\ep_1 = 2W$.
By Lemma \ref{Uniform collar}, for each $t$ there is a
nonperipheral curve $\beta_t$ 
with a collar neighborhood of $q_t$-width $W$.
Since for $t\in J$,  $|\alpha|_{q_t} \le 2W$, we may conclude
that
$i(\alpha,\beta_t)=0$. 
Hence for any $t,s\in J$, $d_\CC(\beta_t,\beta_s) \leq 2$.

The existence of the collar implies
$Ext_{x_t}(\beta_t)\leq 1/W^2$. Hence by Lemma \ref{short are
close},
we conclude $d_\CC(\beta_t,F(t)) \le 1/W^2 + 1$ for $t\in J$.

It follows that $d_\CC(F(s),f(t)) \le 2/W^2 + 4$, and we set
$d_1$
accordingly. 
\end{pf}

%

The following lemma will allow us to convert a long
almost-horizontal
arc which has small diameter (i.e. winds around tightly) to one
which
has a definite diameter, within bounded distance in $\CC$.
\begin{lemma}{horizontal expand}
There exist constants $\ep_3>\ep_2>0$ and $d_2>0$, depending 
on the topology of $S$ and the initial choice of $\ep,\theta$,
so that the following holds. 
Let $\tau$ be an almost horizontal straight segment with respect to $q$,
of
length $|\tau|_q\ge \ep_3$.
Let $J$ be the interval
$$
J = \{t\ge 0: \diam_{q_t}(\tau) < \ep_2\}
$$
and suppose $0\in J$. Then $\diam_\CC(F(J))\le d_2$.
\end{lemma}

\begin{pf}
For any $t$ let $\beta_t$ be the homotopy class of the core of
the
annulus of width $W$ given by Lemma \ref{Uniform collar}. 
Let $\ep_2 = W-2\ep$. 

Let $\Psi^t$ denote the horizontal domain $\Psi_\ep(\tau)$ with
respect to $q_t$ (see Definition \ref{vertical domain def}).
Then $\diam_{q_t}(\Psi^t) \le \diam_{q_t}(\tau) + 2\ep < W$ for
$t\in
J$. Thus any closed curve in $\Psi^t$ has 0 intersection number
with
$\beta_t$. We will show that if $|\tau|_q$ is sufficiently
large,
there exists a nontrivial, nonperipheral curve $\kappa$ which is
contained
in $\Psi^t$ for all $t\in J$.

We use an argument similar to
the proof of Lemma \ref{Vertical domain properties} part
(\ref{non-simple}).  Recall the vertical intervals $\sigma_x$
of radius $\ep$ around
nonsingular $x\in \tau$ from the definition of $\Psi^0$. 
For any $y\in S$ let $f(y) = \#\{x\in \tau: y\in \sigma_x\}$.
Then $\int_{\Psi^0} f(y)$ with respect to  $q$-area is 
$2\ep|\tau|_{q,h}$. 
On the other hand the integral is at most
$\Area_q(\Psi^0)\max f$, so that 
$$
\max f \ge {2\ep\over \Area_q(\Psi^0)} |\tau|_{q,h} \ge
{\ep\over \Area_q(\Psi^0)} |\tau|_q,
$$
where the second inequality is due to  $\tau$ being almost
horizontal.
We also have  $\Area_q(\Psi^0) \le
n_3\diam_q(\Psi^0)^2 \le n_3 W^2$ where $n_3$ was defined
in \S\ref{basic qd}. Thus we have
$ \max f \ge \ep |\tau|_q/ (n_3 W^2) $. 
Set $\ep_3 = 3n_3W^2/\ep$, and now $|\tau|_q\ge\ep_3$
implies
$\max f \ge 3$. As in Lemma \ref{Vertical domain properties} we 
conclude that $\Psi^0$ contains
a nontrivial nonperipheral curve $\kappa$.

Since lengths in the vertical direction shrink as $t$ increases,
for all $t\ge 0$ we have $\Psi^0 \subset \Psi^t$. Thus $\kappa$
is in all the $\Psi^t$, and hence must have 0 intersection number
(hence $\CC$-distance 1) with all $\beta_t$ for $t\in J$, as
above.  As in Lemma \ref{Short Curves Become Balanced},
$d_\CC(\beta_t,F(t))\le 1/W^2 + 1$, so it follows that $F(t)$ is
within bounded $\CC$-distance of $\kappa$ for all $t\in J$.
\end{pf}


The next lemma shows that if an almost vertical straight segment
misses an almost horizontal segment of definite length, then
after a bounded wait as measured in the curve complex, it will
either be very short, or almost horizontal.

\begin{lemma}{Segment Disjoint From Horizontal}
 In addition to our previous
assumptions suppose we also have 
$\ep < \min(\ep_2/a_0, \ep_2a_1/a_2)$. There is  a number
$d_3=d_3(\epsilon,\theta)$ with the following property. Suppose 
 $\alpha$ is a straight segment
of $q$ disjoint from an almost horizontal straight segment $\tau$
of
length
$\epsilon$, let 
\begin{align*}
J &= \{ t\ge 0: \text{$|\alpha|_{q_t}>\ep_1$ and $\alpha$ not
almost
horizontal with respect to $q_t$} \}.
\end{align*}
 Then $\diam_\CC(F(J)) \le d_3$.
\end{lemma}

\bfheading{Remark.} Since $\ep_2$ was given as $W-2\ep$ in Lemma
\ref{horizontal expand}, it is evident that our added conditions
of the form $\ep<C\ep_2$ are satisfied for $\ep$ sufficiently
small. 

\begin{pf}
If $\alpha$ is not almost vertical, then for a bounded $T$
(depending on $\ep,\theta$) it will be almost horizontal with
respect to $q_T$.  By Lemma \ref{Lipschitz}, $F([0,T])$ has
bounded diameter in $\CC$.


Thus we may assume $\alpha$ is almost vertical to begin.  Suppose
its
length is at most $L_0$.  Then for $T=\log 2L_0/\epsilon_1$, its
$q_T$-vertical length is at most $\ep_1/2$; thus it either has
$q_T$-length
less than $\epsilon_1$ or it is not almost vertical.  In the
first
case we have satisfied the conclusion of the Lemma, again
bounding
$\diam_\CC(F([0,T]))$ by Lemma
\ref{Lipschitz}.  In the
second case we are also done by the argument in the first
paragraph.

Thus finally assume $\alpha$ is almost vertical and has length
greater
than $L_0$.  
Since $|\tau|_q = \ep$ and $\tau$ is almost horizontal, for
$t_1=\log 2\ep_3/\ep$ we have  
$|\tau|_{q_{t_1}} \ge \ep_3$. Lemma \ref{horizontal
expand} implies that either $\diam_\CC(F([t_1,\infty)))\le d_2$
in
which case we are done, or
there is a $t_2>t_1$, with
$\diam_\CC(F([t_1,t_2]))\le d_2$, so that
$\diam_{q_{t_2}}(\tau)\ge \ep_2$.


Now if $\alpha$ is not almost vertical or has length at most
$L_0$ 
with
respect to $q_{t_2}$, we are done by the above cases. 
Otherwise, we construct the vertical domain
$\Omega_1=\Omega_\ep(\alpha)$ with respect to $q_{t_2}$. 
Since $\alpha$ is disjoint from $\tau$ and (by assumption) $\ep_2
>
a_0\ep$, Lemma \ref{Vertical domain properties} gives a
subarc $\tau_1$ of $\tau$ of length $a_1\ep_2$,
disjoint from $\Omega_1$.
The total horizontal
length of $\boundary\Omega$ is bounded by 
$a_2\ep$ by part (\ref{horizontal
boundary length}) of Lemma \ref{Vertical domain properties}. 
Thus, since $a_2\ep < a_1\ep_2$ and  applying
Lemma \ref{horizontal segment in disk}, we conclude
that the component $Y$ of $S\setminus\Omega$ containing $\tau_1$
cannot be
a disk or once-punctured disk. 
Thus $\boundary \Omega_1$ has 
nontrivial and nonperipheral components. For each such
component $\sigma$, we have
$|\sigma^*|_{q_{t_2}} \le \ell_0 = a_2\ep + 1/\ep$ by part
(\ref{boundary length}) of Lemma \ref{Vertical domain
properties},
where $\sigma^*$ is the geodesic representative.  
This bound means that within bounded Teichm\"uller distance
either
$|\sigma^*|_{q_t}$ reaches $\ep_1$, or it starts to increase.
Applying
Lemma \ref{Short Curves Become Balanced}, we conclude that either
the
remaining 
$\diam_\CC(F([t_2,\infty)))$ is bounded, in which case we are
done, or
there is $t_3$ with 
bounded $\diam_\CC(F([t_2,t_3]))$, such that for $t>t_3$,
$|\sigma^*|_{q_t} > \ep_1$. It follows that for an additional
$t_4$
with $t_4-t_3$ bounded, $|\sigma^*|_{q_{t_4}} \ge 2\ell_0$.

We can now repeat the argument:
There is a $t_5$ such that $\diam_\CC(F([t_4,t_5]))$ is bounded,
so that
with respect to $q_{t_5}$ we either have the desired condition
for
$\alpha$, or $\alpha$ is almost vertical, of length at least
$L_0$,
and $\tau_1$ now  has diameter at least $\ep_2$.
Thus Lemma \lref{Vertical domain properties} again gives 
a vertical domain $\Omega_2 = \Omega_\ep(\alpha)$ with respect to
$q_{t_5}$ whose
boundary components have $q_{t_5}$-length at most $\ell_0$.
Since our previous boundary components $\sigma$ 
now have $|\sigma^*|_{q_{t_5}}\ge 2\ell_0$,  
no nontrivial nonperipheral component of $\boundary\Omega_2$ is
homotopic to $\sigma$. 
The vertical domains decrease monotonically as $t$ increases, 
so we conclude that $i_*(\pi_1(\Omega_2))$ is a proper subgroup
of
$i_*(\pi_1(\Omega_1))$. 

We may repeat this procedure, obtaining a sequence $\Omega_{j+1}
\subset \Omega_j$ which terminates in at most $n_4$ steps, at
which
point $\alpha$ has length less than
$\epsilon_1$ or is almost horizontal, as desired, or we find that
a
remaining interval $[t,\infty)$ has bounded-diameter image. 
\end{pf}

 From now on let us assume that $\ep,\theta$ satisfy the
conditions of
Lemma \ref{Segment Disjoint From Horizontal} as well as the
previous
conditions. 

The following lemma shows that if a curve is balanced at
$L_q(0)$,
then in the forward direction it will become almost horizontal
after an interval of bounded size in the curve complex.

\begin{lemma+}{Almost Horizontal}  There exists $d_4=
d_4(\ep,\theta)$  so that if  $\beta$ is balanced
with respect to $q$ and  $$ J = \{t\ge 0: \text{$\beta$ is not
almost horizontal with respect to $q_t$}\},  $$ then 
$\diam_\CC(F(J)) \le d_4$.
\end{lemma+}
\begin{pf}  Let $n_3$ be the bound 
for the number of mutually disjoint saddle
connections one can have in $S$. 
Thus, $\beta$ runs through at most $n_3$ saddle connections, 
although some may be traversed arbitrarily many times. 
 
Since $|\beta|_{q,h} = |\beta|_{q,v}$, for $t>t_1 = \half\log
1/\theta$ we have $|\beta|_{q_{t},h} > {1\over\theta}
|\beta|_{q_{t},v}$.  It follows that if $\beta_{0,t}$ is the
subset of $\beta$ which traverses almost-horizontal arcs with
respect to $q_{t}$, we have $|\beta_{0,t}|_{q_{t}} \ge
|\beta\setminus\beta_{0,t}|_{q_{t}}$. Now by Lemma \ref{Short
Curves
Become Balanced}, we have $t_2$ with bounded
$\diam_\CC(F([t_1,t_2]))$ such that $|\beta|_{t_2} \ge \ep_1$,
and therefore $|\beta_{0,t_2}|_{q_{t_2}} \ge \ep_1/2$. It may
still be that this length is obtained by traversing many times a
very short almost horizontal curve $\sigma$, but applying Lemma
\ref{Short Curves Become Balanced} again we obtain $t_3$ with
bounded $\diam_\CC(F([t_2,t_3]))$ such that $|\sigma|_{q_{t_3}}
\ge \ep_1$. 

In particular $\beta$ contains an embedded straight segment
$\tau$ which is
almost-horizontal with respect to $q_{t_3}$ 
and of length at least $\ep$.

Now suppose $\beta$ is not almost horizontal for $q_{t_3}$, so
that there is a segment $\beta_1$ of $\beta$ which has length at
least $\epsilon$ and is not almost horizontal. Since $\beta$ has
no self intersections, $\beta_1$ is disjoint from $\tau$.
Applying Lemma \lref{Segment Disjoint From Horizontal}, there
exists $t_4$
with $\diam_\CC(F([t_3,t_4])) \le d_3$, such that if $t>t_4$ then
either $\beta_1$ is
almost horizontal with respect to $q_{t}$, or
$|\beta_1|_{q_{t}}\leq \epsilon_1$. In the latter case, set
$t_5 = t_4 + \log 2\epsilon_1/\ep\theta$, and note that for
$t>t_5$,
$\beta_1$ will either have length less than $\ep$ or be almost
horizontal.  Apply this to all of the saddle connections of
$\beta$.
\end{pf}

The next lemma shows that, unless a curve is almost vertical in
$L_q(0)$, it can be balanced in the forward direction after an
interval of bounded size in the curve complex. 
\begin{lemma+}{Almost Vertical} There exists
$d_5=d_5(\ep,\theta)$ such that 
if  $\gamma$ is
not almost vertical with respect to $q$ then for
$$
J  = \{ t\ge 0: |\gamma|_{q_t,v} > |\gamma|_{q_t,h}\}
$$
we have $\diam_\CC(F(J)) \le d_5$.
\end{lemma+}

\begin{pf}  
Since $\gamma$ is not almost vertical with respect to $q$, 
it contains a segment $\tau$ that
is not almost vertical and has length at least $\epsilon$. 
For $t_1 = \log 2/\theta$, $\tau$ will be
almost horizontal and have length at least $\epsilon$ with
respect to $q_{t_1}$. 
We can assume that there is a set of almost vertical
saddle connections  $\omega\subset\gamma$ 
that carry at least $1/2$ of the length of $\gamma$, for
otherwise
$|\gamma|_{q_t,h}$ would dominate for $t>t_2$ for a bounded
$t_2$.
Obviously each $\omega$  is disjoint from $\tau$,
since $\gamma$ does not have self intersections. 

By Lemma \ref{Segment Disjoint From Horizontal}, there is
$t_3>t_1$
with $\diam_\CC(F([t_1,t_3]))\le d_3$ such that either
$\omega$ is almost horizontal or has length at most $\epsilon_1$
with
respect to $q_{t_3}$.
If all the $\omega$'s are in the
former case we are done since then $\gamma$ would have been
balanced
for $t\le t_3$.

Suppose, then, that
some $\omega$ has length smaller than $\epsilon_1$.
If we follow any strand of $\gamma$ starting at $\omega$ until it
returns with the same orientation, we may find in the set of
saddle connections it traverses a geodesic loop homotopic to a
simple loop.  There is a bound (in terms of the number of possible saddle
connections $n_2$) on the number of such loops $\gamma'$, and hence
there must be some $\gamma'$ each of whose saddle connections are
traversed at least a definite fraction of the number of times
$\omega$ is traversed.

Either $\gamma'$ contains an
almost horizontal saddle connection, or it has length at most
$n_2\ep_1$. In that case, within bounded $t>t_3$ it will either have
length $\ep_1$ or begin to grow, and we may apply Lemma 
\ref{Short Curves Become Balanced} to get a $t_4$ 
with bounded $\diam_\CC(F([t_3,t_4]))$, 
such that for $t>t_4$
$|\gamma'|_{q_t}\ge \ep_1$ and
$|\gamma'|_{q_t}$ is growing, which implies that
after an additional bounded $t_5$, it will be balanced. 

Thus for $t>t_5$, the 
contribution of $\omega$ to $\gamma$ is either itself horizontal
or offset (to within a bounded factor) by the other segments in
$\gamma'$, and after another bounded interval $[t_5,t_6]$ we have
balance.
\end{pf}

\subsection{Completion of proof}
\label{finishproof}
We must show that all three conditions of definition
\lref{contraction property} hold for our projection $\pi_q$, with
suitable constants $a,b$ and $c$ (independent of the geodesic
$L_q$).

In what follows, fix $q$ and $x = L_q(0)$. Let $\alpha = F_q(0)$ be a
shortest curve on $x$. Assume $\pi_q(\beta) = 0$ so that $\beta$ is
balanced at $0$. We assume that all curves are $q$-geodesics (hence
$q_t$-geodesics for any $t$), and furthermore that $\ep,\theta$ satisfy
the conditions in Lemma \ref{beta and gamma intersect} and Lemma
\ref{Segment Disjoint From Horizontal}.

Let us restate the conditions in our current
terminology, and prove them. 

\bfheading{Condition (1):} $\diam_\CC(F_q([0,\pi_q(F_q(0))])) \le
c$.

We may assume $|\alpha|_{q,v} > |\alpha|_{q,h}$, or equivalently
that
the balance point $\pi_q(\alpha)$ is positive.
Since $\alpha$ has minimal extremal length with respect to $x$,
we have a bound $|\alpha|_q \le Ext_x(\alpha)^{1/2} \le
\sqrt{e_0}$.
Thus for bounded $t_1$, the vertical length
$|\alpha|_{q_{t_1},v}$
becomes $\ep_1$, so either $\alpha$ is balanced for $t\le t_1$,
in
which case we are done, or $|\alpha|_{q_{t_1}}\le \ep_1$. In the
latter case we apply 
Lemma \lref{Short Curves Become Balanced} to see that 
either $\alpha$ is vertical, in which case
$\diam_\CC(F([0,\infty)))$
is bounded and $\pi_q(\alpha) = +\infty$, so we are done,
or there is $t_2$ with $\diam_\CC(F([t_1,t_2]))$ bounded, and
$|\alpha|_t$ is increasing after $t_2$, so it is balanced for some
$t<t_2+\half\cosh^{-1}\sqrt 2$ and again we are done.

\bfheading{Condition (2):}
If $d_\CC(\beta,\gamma) \le 1$ then
$\diam_\CC(F_q([\pi_q(\beta),\pi_q(\gamma)])) \le c$.

Recall $\pi_q(\beta) = 0$. Assume $\pi_q(\gamma)>0$, so 
$|\gamma|_{q,v} > |\gamma|_{q,h}$.
Since $\beta$ is balanced at 0, 
Lemma \ref{Almost Horizontal} gives $t_1>0$ with
$\diam_\CC(F_q([0,t_1]))$ bounded, such that $\beta$ is almost
horizontal with respect to $q_{t_1}$. Lemma \ref{Short Curves
Become
Balanced} then gives $t_2$ with $\diam_\CC(F_q([t_1,t_2]))$
bounded,
so that $|\beta|_{q_{t_2}} \ge \ep_1$. Thus there is a $t_3\ge
t_2$ with 
$F_q[t_2,t_3]$ bounded, so
that $\beta$ has an almost horizontal segment
of length $\ep$ with respect to $q_{t_3}$. 

Now, either $\gamma$ is already balanced for $t\le t_3$, in which
case we
are done, or it is still mostly vertical with respect to
$q_{t_3}$. In
this case, since $\gamma$ is disjoint from $\beta$ it misses the
horizontal segment of length $\ep$ and Lemma \ref{Segment
Disjoint From Horizontal} gives $t_4$ with $\diam_\CC(F_q([t_3,t_4]))$
bounded,
so that every saddle connection of $\gamma$ is either almost
horizontal or has length at most $\ep_1$ with respect to
$q_{t_4}$.
Thus for $t_5$ with bounded $t_5-t_4$, either the length of
$\gamma$
shrinks to $\ep_1$, or it begins to increase so $\gamma$ is
balanced
for $t\le t_5$. In the former case, Lemma \ref{Short Curves
Become
Balanced} says that the segment $J$ of $t>t_5$ where
$|\gamma|_{q_t}\le \ep_1$ has bounded-diameter image in $\CC$.
This
includes the case where $\gamma$ is completely vertical and
$\pi_q(\gamma) = +\infty$. In all other cases, $\gamma$ will be
balanced for some $t\in J$ and again we are done. 

\bfheading{Condition (3):}
If $d_\CC(\beta,F_q(\pi_q(\beta))) \ge a$ and
$d_\CC(\beta,\gamma) \le
b d_\CC(\beta,F_q(\pi_q(\beta)))$ then
$ \diam_\CC F_q([\pi_q(\beta),\pi_q(\gamma)]) \le c$.

Recall that $F_q(\pi_q(\beta)) = F_q(0) = \alpha$.
Assume without loss of generality that $\gamma$ is more vertical
than
horizontal at $q_0$. 
By Proposition \ref{Distance 1}, we can assume that
$\beta\in \CC_b$, for otherwise its 
distance from the image of $F$ is at most $1$.  By Lemma
\lref{Almost
Horizontal}, there is some $t_1>0$ with $\diam_\CC(F_q[0,t_1])$ 
bounded such that $\beta$ is almost horizontal 
at $q_{t_1}$. 

We next show that $\gamma$ cannot be almost vertical at
$q_{t_1}$.  Let $M,h$ be the constants given by lemma \ref{beta
and gamma intersect}, and suppose by contradiction that $\gamma$
is almost vertical.  If $\alpha'$ is any curve of
$q_{t_1}$-length at most $M$, there is a bound $d(M)$ on
$d_\CC(\alpha',\alpha)$ by the following: Lemma \ref{Uniform
collar} gives a nonperipheral annulus $A$ of width $W$ so that
the intersection of $\alpha'$ with its core $\sigma$ is at most
$M/W$. Lemma \ref{connected} then bounds
$d_\CC(\alpha',\sigma)$. Lemma \ref{short are close} in turn
bounds $d_\CC(\sigma,F_q(t_1))$ since both have bounded extremal
length.  Finally $d_\CC(F_q(t_1),F_q(0)) =
d_\CC(F_q(t_1),\alpha)$ is bounded by choice of $t_1$.

Now applying Lemma \ref{nesting consequence} with $Q=2/h$, and 
$k = d(M)$, we obtain, provided $d_\CC(\alpha,\beta) \ge D_3$ and
$d_\CC(\gamma,\beta) \le \nu d_\CC(\alpha,\beta)$, that $$
\min_{\alpha'} i(\beta,\alpha') \min_{\alpha'} i(\gamma,\alpha') 
\ge Q i(\beta,\gamma)
$$
where $\alpha'$ varies over all curves of $q_{t_1}$-length at
most
$M$.
On the other hand, Lemma \ref{beta and gamma intersect} gives the
opposite inequality
$$
i(\gamma,\beta) \ge h \min_{\alpha'} i(\beta,\alpha')
\min_{\alpha'}
i(\gamma,\alpha'). 
$$
This is a contradiction since we have chosen $Q>1/h$,
and we conclude that $\gamma$ cannot be almost vertical at $t_1$.
Thus by Lemma \lref{Almost Vertical}, there exists $t_2$ with 
$\diam_\CC(F_q([t_1,t_2]))$ bounded such that $\gamma$ is
balanced at
$t_2$. This concludes the proof of Theorem \ref{Projection Theorem}.

\section{Contraction property and hyperbolicity}
\label{hyperbolic}

To complete the proof of Theorem \ref{Hyperbolicity}, it remains to
prove Theorem \ref{contraction implies hyperbolicity}, that if
a geodesic metric space 
$X$ has a coarsely transitive path family $\Gamma$ with the
contraction property then $X$ is hyperbolic.

For our purposes a path $\gamma:I\to X$ is a {\em quasi-geodesic} if
the following inequality holds for any $x,y\in I$:
$$
\length_s(\gamma[x,y]) \le 
K d_X(\gamma(x),\gamma(y)) + \delta
$$
where $K\ge 1$ and $\delta,s\ge 0$ are fixed constants, and
$\length_s$ for $s>0$ is
``arclength on the scale $s$'', which is defined as follows:
$\length_s(\gamma[x,y]) = sn$ where $n$ is the smallest number for
which $[x,y]$ can be subdivided into $n$ closed subintervals
$J_1,\ldots,J_n$ with $\diam_X(\gamma(J_i)) \le s$.
(This definition circumvents the need for checking the behavior of  
the parametrization at small scale; we let $\length_0$ denote normal length).
Note also that the opposite inequality $d_X(\gamma(x),\gamma(y)) \le
\length_s(\gamma[x,y])$ holds automatically.

The proof is in two steps.
We say that $X$ has {\em stability of quasi-geodesics} if for all
$K\ge1$,$\delta,s\ge 0$ there exists $R>0$ such that any
$(K,\delta,s)$-quasi-geodesic 
$\alpha:I\to X$ with endpoints $x,y$ remains in an $R$-neighborhood
of any geodesic $[xy]$. 

\begin{lemma}{contraction implies stability}
If $X$ has a coarsely transitive path family $\Gamma$ with the contraction
property then $X$ has stability of quasi-geodesics. In addition, the
paths of $\Gamma$ themselves are uniform quasi-geodesics.
\end{lemma}

\begin{lemma}{stability implies hyperbolicity}
Stability of quasi-geodesics implies hyperbolicity.
\end{lemma}

\begin{pf*}{Proof of Lemma \ref{contraction implies stability}}
We may assume that the path family $\Gamma$ is transitive, since 
for paths of length bounded by a fixed $D$ it is easy to define a
contracting projection, simply by mapping all of $X$ to one endpoint. 

Consider $\gamma:[0,M]\to X$ in $\Gamma$, and let $\alpha:[0,L]\to X$
be a $(K,\delta,s)$-quasi-geodesic such that $\alpha(0) = \gamma(0)$
and $\alpha(L) = \gamma(M)$. Note that a $(K,\delta,0)$-quasi-geodesic
is also a $(K,\delta+s,s)$-quasi-geodesic for any $s>0$ since
$\length_s \le \length_0 + s$. Thus from now on we may assume $s>0$.

We show that $\alpha$ remains in a
$R(K,\delta,s)$-neighborhood of $\gamma$.  
The proof is somewhat
complicated by the fact that we do not assume continuity of $\alpha$,
$\gamma$ or $\pi$, but the idea is simple and well-known: large
excursions of $\alpha$ away from $\gamma$ can, using the contraction
property, be circumvented by short cuts that travel along the
projection to $\gamma$.

Let $r(u) = d(\alpha(u),\gamma(\pi(\alpha(u))))$. We will bound $r(u)$
uniformly in terms of $K,\delta,s,$ and the constants $a,b$ and $c$ of
the contraction property (Definition \ref{contraction property}).

Divide $[0,L]$ into closed intervals $J_1,\ldots,J_n$ such that 
$ns = \length_s(\alpha[0,L])$ and $\diam \alpha(J_i) \le s$. Then by
part (\ref{quasi-lipschitz}) of definition \ref{contraction  
property},
$\diam \gamma(\pi(\alpha(J_i))) \le s'$, where $s' = c$ if $s\le 1$
and $s' = 1+cs$ if $s>1$. 

Fix $R_0>0$, to be determined shortly.
For any $u\in [0,L]$, 
if $r(u) \ge R_0 + s'$ then $u$ is contained in some interval  
$J=[u_0,u_1]$, a
union of $J_i$, such that $r\ge R_0$ in $J$ and $r\le R_0+s'$ at  
$u_0$
and $u_1$.
Subdivide $J$ into intervals $K_1,\ldots,K_m$, each a union of at  
most
$bR_0/s$ of the $J_i$,
so that for each $j$ $\diam \alpha(K_j) \le bR_0$, and the number
$m$ is at most $1+\length_s(\alpha(J))/bR_0$. Now assuming $R_0 \ge  
a$
and applying the
contraction property (part \ref{contraction}) to each of these we
obtain
$$
\diam \gamma([\pi(\alpha(u_0)),\pi(\alpha(u_1))])\le
mc
$$
and by the triangle inequality
\begin{equation}
\label{distance bound}
d(\alpha(u_0),\alpha(u_1)) \le 
2(R_0+s') + 
\left(1+{\length_s(\alpha(J))\over bR_0}\right)c.
\end{equation}
Since $\alpha$ is a $(K,\delta,s)$-quasi-geodesic,  
$\length_s(\alpha(J))
\le Kd(\alpha(u_0),\alpha(u_1))+\delta$. Combining with  
(\ref{distance
bound}), we get 
\begin{equation}
\label{length bound}
\length_s(\alpha(J)) \le {Kc\over bR_0}\length_s(\alpha(J)) + 
2Kc(R_0+s') + \delta.
\end{equation}
Make the (a priori) choice of 
 $R_0$ sufficiently large that $Kc/bR_0< 1/2$. Then (\ref{length
bound}) gives an upper bound $R$ on $\length_s(\alpha(J))$ depending
only on the initial constants.

Thus $d(\alpha(u),\{\alpha(u_0),\alpha(u_1)\})$ is at most $R/2$,
and in particular, no point in $\alpha(J)$ can be further
than $R_0 + R/2$ from $\gamma([0,M])$. 
Furthermore by applying part (\ref{quasi-lipschitz}) of the
contraction property, it follows that $r(u)$ is bounded uniformly.

This implies that we can project from $\gamma$ back to $\alpha$, in
the following sense:
For any $t\in[0,M]$ we can find $u\in[0,L]$ such that
$d(\gamma(t),\gamma(\pi(\alpha(u))))$ is bounded by a uniform
constant, just by chopping $\alpha$ into bounded-length pieces
and applying parts (1) and (2) of the contraction property.
Now by the bound on $r(u)$ we can bound $d(\gamma(t),\alpha(u))$
uniformly.

Apply this to an actual geodesic $\alpha$ and a quasi-geodesic  
$\beta$
with the same endpoints. Letting $\gamma\in \Gamma$ be a path with  
the
same endpoints, we project from $\beta$ to $\gamma$ and then from
$\gamma$ to $\alpha$ as above. Both steps move a bounded distance,  
so
we conclude that $\beta$ lies in a
bounded neighborhood of $\alpha$. Hence, we have stability of  
quasi-geodesics.
\end{pf*}

\begin{pf*}{Proof of Lemma \ref{stability implies hyperbolicity}}
To prove hyperbolicity it suffices to establish the thin triangle
condition. Let $x,y,z$ be three points in $X$. We must show that
$[xy]$ lies in a $\delta$-neighborhood of $[xz]\union[yz]$, for
uniform $\delta$.

Let $z'\in [xy]$ be a point that minimizes distance from $z$ to  
$[xy]$.
We claim that the broken geodesic $[xz']\union[z'z]$ is a
$(3,0,0)$-quasi-geodesic. 
If $z' = x$ this is obvious, so assume $z' \ne x$.
Let $u$ lie in $[xz']$ and $v$ lie in $[z'z]$. 

It follows from the choice of $z'$ that it also minimizes 
distance from $[v]$ to $[xy]$ (via the triangle inequality).
Thus $d(u,v) \ge d(z',v)$.

By the triangle inequality, $d(u,v) \ge d(u,z') - d(z',v)$. Thus
adding this to twice the previous inequality we get
$3 d(u,v) \ge d(z',v) + d(u,z')$. 
This is exactly the fact that $\length([uz']\union[z'v])$ estimates
$d(u,v)$, so we conclude $[xz']\union [z'z]$ is a  
$(3,0,0)$-quasi-geodesic.

Now by stability of quasi-geodesics, we have  that
$[xz']\union [z'z] $ is in a uniform $\delta$-neighborhood of $[xz]$,  
and in
particular $[xz']$ is. Applying the same argument for $y$ replacing
$x$, we see that all of $[xy]$ is in a $\delta$-neighborhood of
$[xz]\union[yz]$. This concludes the proof.
\end{pf*}

\section{Relative Hyperbolicity}
\label{relative hyperbolicity}
In this final section we establish Theorems \ref{Relative Hyperbolicity 1}
and \ref{Relative Hyperbolicity 2}, which provide an interpretation of our
hyperbolicity theorem in terms of the geometry of Teichm\"uller
space, and the structure of the Mapping Class Group.

The following terminology is due to Farb \cite{farb:thesis}: If
$X$ is any geodesic metric space and $\HH$ is a family of regions
in $X$, let the {\em electric distance} $d_e$ on $X$ be the path
metric imposed by shrinking each $H\in\HH$ to diameter 1, in the
following way: For each $H\in\HH$ create a new point $c_H$ and an
interval of length $1/2$ from $c_H$ to every point in $H$. The
new metric is induced by shortest paths in this enlarged space
$\hat X$ (called the {\em electric space}).  We say $X$ is {\em
relatively hyperbolic} with respect to $\HH$ if $(\hat X,d_e)$ is
$\delta$-hyperbolic for some $\delta$.

\subsection{In Teichm\"uller space}
Fixing $\ep_0>0$ sufficiently small that the Collar Lemma holds
for $\ep_0$, let $\HH_C = \{H_\alpha\}_{\alpha\in \CC_0(S)}$
denote the family of regions in $\TT(S)$ defined as in the
introduction:
$$
H_\alpha = \{x\in\TT(S): Ext_x(\alpha) < \ep_0\}.
$$
Then it is easy to see that a set of points
$\alpha_1,\ldots,\alpha_k$ is a simplex in $\CC(S)$ if and only
if $H_{\alpha_1}\intersect\cdots\intersect H_{\alpha_k}$ is
non-empty. In other words, $\CC(S)$ is the {\em nerve} of the
family $\HH_C$.

The statement of Theorem \ref{Relative Hyperbolicity 1} is a direct
consequence of Theorem \ref{Hyperbolicity} and the following:

\begin{lemma}{nerve is quasi}
The electric space $(\hat\TT(S),d_e)$ defined with respect to
the family $\HH_C$ is quasi-isometric to $\CC_1(S)$.
\end{lemma}

\begin{pf}
There is a natural map $\varphi: \CC_0(S)\to \hat\TT(S)$ taking
each $\alpha$ to the new point $c_\alpha\equiv c_{H_\alpha}$.  The set
$\CC_0(S)$ is clearly $1/2$-dense in $\CC_1(S)$. Let us check
that its image $\{c_\alpha\}$ is $d_0$-dense in $\hat\TT(S)$,
for some $d_0<\infty$.

Recall that for any conformal structure
$x$ on $S$ there is a curve $\alpha\in\CC_0(S)$ with $Ext_x(\alpha)
\le e_0$. Then for this $\alpha$, we see that $x$ is a bounded
Teichm\"uller distance
(in fact $\half\log(e_0/\ep_0)$) from $H_\alpha$: 
we may apply to $x$ a Teichm\"uller map whose vertical
foliation consists of leaves homotopic to $\alpha$. It follows
that $\{c_\alpha\}$ is $(\half + \half\log(e_0/\ep_0))$-dense
in $\hat\TT(S)$.

Now we need only show that for any $\alpha,\beta\in\CC_0(S)$
\begin{equation}
\label{quasi isom def}
{1\over K}d_\CC(\alpha,\beta) - a \le d_e(c_\alpha,c_\beta)
			\le Kd_\CC(\alpha,\beta) + a
\end{equation}
with fixed $K,a>0$, to show that $\varphi$ induces a
quasi-isometry.  One direction is easy: if $d_\CC(\alpha,\beta)
= 1$ then $H_\alpha\intersect H_\beta$ is nonempty, and any point
$x$ in this set is connected to each of $c_\alpha$ and $c_\beta$
by a segment of length $1/2$. Hence $\varphi$ is $1$-Lipschitz.

To obtain the other direction, 
consider for any $x\in\TT(S)$ the set $\Phi(x)$ of elements in
$\CC_0(S)$ of minimal $Ext_x$.
This set has diameter at most $2e_0+1$ by Lemma \ref{Phi diam bound}.
Now if $d_\TT(S)(x,y) \le 1$ we see also that $\Phi(x)\union \Phi(y)$
has bounded diameter, by Lemma \ref{Lipschitz}.
Note also that if $x\in H_\alpha$ then $d_\CC(\alpha,\Phi(x))\le 1$.

Thus, any map that associates to $x\in\TT(S)$ some (any) element
of $\Phi(x)$ and to $c_\alpha$ associates $\alpha$ will expand distances
by a bounded multiplicative and additive amount, and serve as an
inverse to $\varphi$. It follows that $\varphi$ is a quasi-isometry.
\end{pf}

\subsection{In the Mapping Class Group}
To carry out a similar analysis for $\Mod(S)$, recall first that
for any group $G$ with a fixed finite generating set $\Gamma$,
the Cayley graph $\GG = \GG_{G,\Gamma}$ is a 1-complex whose
vertex set is $G$ and whose edges are all pairs $(g,g\gamma)$
with $\gamma\in\Gamma$. Giving all edges length 1, we obtain a
complete locally finite geodesic metric space.

Now for $G=\Mod(S)$, we single out a number of subgroups as
follows.  Up to the action of $\Mod(S)$, there are only a finite
number of distinct non-trivial non-peripheral homotopy classes of
simple curves in $S$ (distinguished by the topological type of
their complement).  Let $\{\alpha_1,\ldots,\alpha_N\}$ be a fixed
list of representatives of these $\Mod(S)$-orbits. Let
$Fix(\alpha_j)$ be the subgroup of $\Mod(S)$ fixing $\alpha_j$.

Given any $\beta\in\CC_0(S)$, let $\alpha_j$ be the unique
representative of $\beta$ in the list, and let $G_\beta$ be
the left-coset of $Fix(\alpha_j)$ defined by $G_\beta =
\{g\in\Mod(S): g(\alpha_j) = \beta\}$.

Now we may form the electric space $\hat\GG$ of $\GG$ relative to
the family of cosets $\{G_\beta\}$, and its electric distance
$d_e$.  The analogue to Lemma \ref{nerve is quasi} is:

\begin{lemma}{coset graph is quasi}
Fixing a choice of generating set $\Gamma$ and representatives
$\{\alpha_1,\ldots,\alpha_N\}$ of $\Mod(S)$-orbits in $\CC_0(S)$,
the electric space $(\hat\GG,d_e)$ is quasi-isometric to
$\CC_1(S)$.
\end{lemma}

Again, this together with Theorem \ref{Hyperbolicity} proves Theorem
\ref{Relative Hyperbolicity 2}, where the relative hyperbolicity of
$\Mod(S)$ is with respect to this family of cosets $\{G_\beta\}$.

\begin{pf}
Let $c_\beta$ denote the new point added to $G_\beta$ in the
construction of $\hat\GG$. The natural map $\varphi:\CC_0(S) \to
\hat\GG$ is again $\varphi(\beta)  = c_\beta$. In this case it is
clear that $\{c_\beta\}$ is $1/2$-dense in $\hat\GG$ since every
$g\in\Mod(S)$ is in the coset $gFix(\alpha_j) = G_{g(\alpha_j)}$ for
each $j\le N$. It remains to check that the inequalities
(\ref{quasi isom def}) hold.

Up to the action of $\Mod(S)$ there are only finitely many pairs
$(\beta,\beta')$ of disjoint curves in $\CC_0(S)$ (i.e. edges in
$\CC_1(S)$).  Let $\{(\beta_i,\beta'_i)\}_{i=1}^L$ be an
enumeration of representatives of $\Mod(S)$-orbits.
For each  $\alpha_j$ there is some (in fact several) $\beta_i$
equivalent to it under $\Mod(S)$, so let $w_{ij}$ be a fixed
group element such that $w_{ij}(\alpha_j) = \beta_i$. Define
$w'_{ij}$ similarly. Since this is a finite list, there is some
upper bound $B$ on their lengths as words in the generating set
$\Gamma$.

Now let $\beta,\beta'\in\CC_0(S)$ be any two curves of distance
1. Hence there exists $g\in\Mod(S)$ and $i\le L$ such that
$g(\beta_i) = \beta$ and $g(\beta'_i) = \beta'$. There also exist
$j,k\le N$ such that $w_{ij}(\alpha_j) = \beta_i$ and
$w'_{ik}(\beta_k) = \beta'_i$. 

Thus $gw_{ij}\in G_\beta$ and $gw'_{ik}\in G_{\beta'}$, and these
two elements are separated by a path in $\GG$ of distance at most
$2B$. We conclude that $d_e(c_\beta,c_{\beta'}) \le 2B+1$ and hence
the map $\varphi$ is $(2B+1)$-Lipschitz.

To obtain a bound in the other direction, note that for any $g\in
\Mod(S)$ we may associate the set $A_g = \{g(\alpha_j)\}_{j\le N}$ in
$\CC_0(S)$, and that the diameter of this set in $\CC(S)$ is
equal to the diameter of $A_{id}=\{\alpha_j\}_{j\le N}$, which is some
fixed $D$ (with appropriate choice of $\alpha_J$ we can easily get $D=2$).
Now given $g$ and $g\gamma$ where $\gamma\in\Gamma$ is
a generator, the distance between the sets $A_g$ and
$A_{g\gamma}$ is equal to that between $A_{id}$ and $A_{\gamma}$,
which is again bounded. Note finally that if $g\in G_\beta$ then
$\beta\in A_g$. Thus we can map the vertices of $\hat\GG$ back to
$\CC_0$, taking each $c_\beta$ to $\beta$, and each $g$ to some (any)
element of $A_g$, and the resulting map is Lipschitz, and inverts $\varphi$.
This proves that $\varphi$ is a quasi-isometry.
\end{pf}


\begin{thebibliography}{10}

\bibitem{ahlfors:invariants}
L.~Ahlfors, {\em Conformal invariants: topics in geometric function theory},
  McGraw-Hill, 1973.

\bibitem{short:notes}
Alonso, Brady, Cooper, Ferlini, Lustig, Mihalik, Shapiro, and Short, {\em Notes
  on word hyperbolic groups}, Group Theory from a Geometrical Viewpoint, ICTP
  Trieste 1990 (E.~Ghys, A.~Haefliger, and A.~Verjovsky, eds.), World
  Scientific, 1991, pp.~3--63.

\bibitem{ballmann:spaces}
W.~Ballmann, {\em Lectures on spaces of nonpositive curvature}, Birkh\"auser,
  1995.

\bibitem{bers:pseudoanosov}
L.~Bers, {\em An extremal problem for quasiconformal mappings and a theorem by
  {T}hurston}, Acta Math. {\bf 141} (1978), 73--98.

\bibitem{bowditch:hyperbolicity}
B.~Bowditch, {\em Notes on {Gromov's} hyperbolicity criterion for path-metric
  spaces}, Group theory from a geometrical viewpoint (Trieste, 1990), World
  Scientific Publishing, 1991, pp.~64--167.

\bibitem{bridson:simplicial}
M.~R. Bridson, {\em Geodesics and curvature in metric simplicial complexes},
  Group Theory from a Geometrical Viewpoint, ICTP Trieste 1990 (E.~Ghys,
  A.~Haefliger, and A.~Verjovsky, eds.), World Scientific, 1991, pp.~373--463.

\bibitem{buser:surfaces}
P.~Buser, {\em {Geometry and Spectra of Compact Riemann Surfaces}},
  Birkh\"auser, 1992.

\bibitem{cannon:negative}
J.~Cannon, {\em The theory of negatively curved spaces and groups}, Ergodic
  theory, symbolic dynamics, and hyperbolic spaces (Trieste, 1989), Oxford
  Univ. Press, 1991, pp.~315--369.

\bibitem{casson:unpub}
A.~J. Casson, {\em Automorphisms of surfaces after {N}ielsen and {T}hurston},
  Notes by S. A. Bleiler, U. T. Austin, 1982.

\bibitem{c-d-p}
M.~Coornaert, T.~Delzant, and A.~Papadopoulos, {\em G\'eom\'etrie et theorie de
  groupes: les groups hyperboliques de {G}romov}, Springer-Verlag, 1990.

\bibitem{farb:thesis}
B.~Farb, {\em Relatively hyperbolic and automatic groups with applications to
  negatively curved manifolds}, Ph.D. thesis, Princeton University, 1994.

\bibitem{travaux}
A.~Fathi, F.~Laudenbach, and V.~Poenaru, {\em Travaux de {T}hurston sur les
  surfaces}, vol. 66-67, Asterisque, 1979.

\bibitem{gardiner}
F.~Gardiner, {\em {T}eichm\"{u}ller theory and quadratic differentials}, Wiley
  Interscience, 1987.

\bibitem{ghys-harpe}
E.~Ghys and P.~de~la Harpe, {\em Sur les groupes hyperboliques d'apr\'es
  {Mikhael Gromov}}, Birkh\"auser, 1990.

\bibitem{gromov:hypgroups}
M.~Gromov, {\em Hyperbolic groups}, Essays in Group Theory {\rm (S. M. Gersten,
  editor)}, MSRI Publications no. 8, Springer-Verlag, 1987.

\bibitem{harer:stability}
J.~Harer, {\em Stability of the homology of the mapping class group of an
  orientable surface}, Ann. of Math. {\bf 121} (1985), 215--249.

\bibitem{harer:cohomdim}
\bysame, {\em The virtual cohomological dimension of the mapping class group of
  an orientable surface}, Invent. Math. {\bf 84} (1986), 157--176.

\bibitem{harvey:boundary}
W.~J. Harvey, {\em Boundary structure of the modular group}, Riemann Surfaces
  and Related Topics: Proceedings of the 1978 Stony Brook Conference (I.~Kra
  and B.~Maskit, eds.), Ann. of Math. Stud. 97, Princeton, 1981.

\bibitem{hatcher}
A.~E. Hatcher, {\em Measured lamination spaces for surfaces, from the
  topological viewpoint}, Topology Appl. {\bf 30} (1988), 63--88.

\bibitem{ivanov:complexes2}
N.~V. Ivanov, {\em Automorphisms of complexes of curves and of {Teichm\"uller}
  spaces}, Preprint.

\bibitem{ivanov:complexes1}
\bysame, {\em Complexes of curves and the {Teichm\"uller} modular group},
  Uspekhi Mat. Nauk {\bf 42} (1987), 55--107.

\bibitem{ivanov:complexes3}
\bysame, {\em Complexes of curves and {Teichm\"uller} spaces}, Math. Notes {\bf
  49} (1991), 479--484.

\bibitem{keen:collar}
L.~Keen, {\em Collars on {R}iemann surfaces}, Discontinuous groups and
  {R}iemann surfaces (Proc. Conf., Univ. Maryland 1973), Ann. of Math. Studies
  79, Princeton, 1974, pp.~263--268.

\bibitem{kerckhoff}
S.~Kerckhoff, {\em The asymptotic geometry of {T}eichm\"uller space}, Topology
  {\bf 19} (1980), 23--41.

\bibitem{kerckhoff:nielsen}
\bysame, {\em The {N}ielsen realization problem}, Ann. of Math. {\bf 117}
  (1983), 235--265.

\bibitem{kerckhoff:simplicial}
\bysame, {\em Simplicial systems for interval exchange maps and measured
  foliations}, Ergodic Theory and Dynamical Systems {\bf 5} (1985), 257--271.

\bibitem{masur:teichgeo}
H.~A. Masur, {\em On a class of geodesics in {T}eichm\"uller space}, Ann. of
  Math. {\bf 102} (1975), 205--221.

\bibitem{masur:wpmetric}
\bysame, {\em The extension of the {Weil-Petersson} metric to the boundary of
  {Teichm\"uller} space}, Duke Math. J. {\bf 43} (1976), 623--635.

\bibitem{masur:iet}
\bysame, {\em Interval exchange transformations and measured foliations}, Ann.
  of Math. {\bf 115} (1982), 169--200.

\bibitem{masur-minsky:complex2}
H.~A. Masur and Y.~Minsky, {\em Geometry of the complex of curves {II}}, in
  preparation.

\bibitem{masur-wolf}
H.~A. Masur and M.~Wolf, {\em Teichm\"uller space is not {G}romov hyperbolic},
  MSRI preprint No. 011-94, 1994.

\bibitem{minsky:2d}
Y.~Minsky, {\em Harmonic maps, length and energy in {T}eichm\"uller space}, J.
  of Diff. Geom. {\bf 35} (1992), 151--217.

\bibitem{minsky:slowmaps}
\bysame, {\em Teichm\"uller geodesics and ends of hyperbolic 3-manifolds},
  Topology {\bf 32} (1993), 625--647.

\bibitem{minsky:extremal}
\bysame, {\em Extremal length estimates and product regions in {T}eichm\"uller
  space}, Duke Math J. {\bf 83} (1996), 249--286.

\bibitem{minsky:taniguchi}
\bysame, {\em A geometric approach to the complex of curves}, Proceedings of
  the 37th Taniguchi Symposium on Topology and {Teichm\"uller} Spaces
  (S.~Kojima et. al., ed.), World Scientific, 1996, pp.~149--158.

\bibitem{minsky:projections}
\bysame, {\em Quasi-projections in {T}eichm\"uller space}, J. Reine Angew.
  Math. {\bf 473} (1996), 121--136.

\bibitem{penner:dilatation}
R.~Penner, {\em Bounds on least dilatations}, Proc. Amer. Math. Soc. {\bf 113}
  (1991), 443--450.

\bibitem{penner-harer}
R.~Penner and J.~Harer, {\em Combinatorics of train tracks}, Annals of Math.
  Studies no. 125, Princeton University Press, 1992.

\bibitem{strebel}
K.~Strebel, {\em Quadratic differentials}, Springer-Verlag, 1984.

\bibitem{wpt:surfaces}
W.~Thurston, {\em On the geometry and dynamics of diffeomorphisms of surfaces},
  Bull. Amer. Math. Soc. {\bf 19} (1988), 417--431.

\bibitem{wolpert:nielsen}
S.~A. Wolpert, {\em Geodesic length functions and the {N}ielsen problem}, J.
  Differential Geom. {\bf 25} (1987), 275--296.

\bibitem{wolpert:plumbing}
\bysame, {\em The hyperbolic metric and the geometry of the universal curve},
  J. Differential Geom. {\bf 31} (1990), 417--472.

\end{thebibliography}
\ifx\undefined\bysame
\newcommand{\bysame}{\leavevmode\hbox to3em{\hrulefill}\,}
\fi

\end{document}